\def\loadTIKZ{\usepackage{tikz}\usetikzlibrary{matrix,arrows,calc,cd,decorations.pathmorphing}}\ifdefined\headpresent\else \providecommand\DocumentClassOptions{article}
\documentclass[\DocumentClassOptions]{amsart}\usepackage{ifpdf}\ifdefined\BibLaTeX \usepackage[\BibLaTeX]{biblatex}\DeclareFieldFormat{labelalpha}{\thefield{entrykey}}\makeatletter\protected\def\abx@missing#1{\mbox{\reset@font\color{red}#1??}}\makeatother \DeclareFieldFormat{title}{\textit{#1}}\DeclareFieldFormat{year}{\textbf{#1}}\renewbibmacro*{date}{\printtext[year]{\printdate}}\BibLaTeXCommands \else\fi \overfullrule=5pt \newcommand\hfuzzReset{\hfuzz=3pt}\hfuzzReset \newcommand\toleranceReset{\tolerance=1400}\toleranceReset \newcommand\emergencystretchReset{\emergencystretch=1ex}\emergencystretchReset \hbadness=10000 \usepackage{xifthen}\usepackage{forarray}\usepackage{xstring}\usepackage{stringstrings}\def\StackCreate#1#2#3{\expandafter\def\csname#1\endcsname{#3}\expandafter\def\csname#1Push\endcsname##1{\expandafter\edef\csname#1\endcsname{##1#2\csname#1\endcsname}}\expandafter\def\csname TopAux#1\endcsname ##1#2##2#3{##1}\expandafter\def\csname#1Top\endcsname{\expandafter\expandafter\expandafter\expandafter\expandafter\expandafter\csname TopAux#1\endcsname\csname#1\endcsname}\expandafter\def\csname PopAux#1\endcsname ##1#2##2#3##3#2{\expandafter\def\csname##3\endcsname{##2#3}}\expandafter\def\csname#1Pop\endcsname{\expandafter\expandafter\expandafter\expandafter\expandafter\expandafter\csname PopAux#1\endcsname\csname#1\endcsname#1#2}}\def\GetAfterColonAux#1:#2;{#2}\def\GetAfterColon#1{\IfBeginWith{#1}{:}{\GetAfterColonAux#1;}{#1}}\usepackage{aliascnt}\newlength{\tempwidth}\newcommand{\fillX}[2][]{\settowidth{\tempwidth}{#2}\def\temp{#1}\ifx\temp\empty\else\addtolength{\tempwidth}{#1}\fi\leavevmode\cleaders\hbox to \tempwidth{\hss #2\hss }\hfill\kern0pt }\usepackage[shortlabels,inline]{enumitem}\setenumerate[1]{leftmargin=5.5ex}\setitemize[1]{leftmargin=5.5ex}\SetEnumitemKey{noindent}{leftmargin=0ex, itemindent=5ex, align=right, itemsep=1ex }\newcounter{ManualLabel}\makeatletter \newcommand\itemPatch[1][]{\item[\theenumi#1]\refstepcounter{ManualLabel}\def\@currentlabel{\theenumi#1}}\newcommand\itemDescribe[1][]{\item[#1]\refstepcounter{ManualLabel}\def\@currentlabel{#1}}\makeatother \ifdefined\selectPages \usepackage[\selectPages]{pagesel} \fi \ifdefined\presentationON \usepackage[hcentering=true,textheight=240mm,textwidth=210mm]{geometry} \else\fi \usepackage{everypage-1x}\ifdefined\AddEverypageHook \newcommand\AddPrivateToMargin[1]{\AddEverypageHook{\tikz[overlay,remember picture]{\node at ($(current page.west)+(1.5,0)$) [rotate=90] {\textcolor{orange}{\vbox{\hrule width \the\textwidth height 0.5pt} \textcolor{defaultcolor}{#1}\ \vbox{\hrule width 40em height 0.5pt}}}; }}}\newcommand\AddLongversionToMargin[1]{\AddEverypageHook{\tikz[overlay,remember picture]{\node at ($(current page.west)+(2,0)$) [rotate=90] {\textcolor{\LongColor}{\vbox{\hrule width \the\textwidth height 0.5pt} \textcolor{defaultcolor}{#1}\ \vbox{\hrule width 40em height 0.5pt}}}; }}}\newcommand\AddOldversionToMargin[1]{\AddEverypageHook{\tikz[overlay,remember picture]{\node at ($(current page.west)+(2.5,0)$) [rotate=90] {\textcolor{\OldColor}{\vbox{\hrule width \the\textwidth height 0.5pt} \textcolor{defaultcolor}{#1}\ \vbox{\hrule width 40em height 0.5pt}}}; }}}\newcommand\AddLineToMargin[3]{\AddEverypageHook{\tikz[overlay,remember picture]{\node at ($(current page.west)+(#2,0)$) [rotate=90] {\textcolor{#1}{\vbox{\hrule width \the\textwidth height 0.5pt} \textcolor{defaultcolor}{#3}\ \vbox{\hrule width 40em height 0.5pt}}}; }}}\else\fi \providecommand\AddPrivateToMargin[1]{\AddToHook {shipout/background}{\tikz[overlay,remember picture]{\node at ($(current page.west)+(1.5,0)$) [rotate=90] {\textcolor{orange}{\vbox{\hrule width \the\textwidth height 0.5pt} \textcolor{black}{#1}\ \vbox{\hrule width 40em height 0.5pt}}}; }}}\providecommand\AddLongversionToMargin[1]{\AddToHook {shipout/background}{\tikz[overlay,remember picture]{\node at ($(current page.west)+(2,0)$) [rotate=90] {\textcolor{\LongColor}{\vbox{\hrule width \the\textwidth height 0.5pt} \textcolor{black}{#1}\ \vbox{\hrule width 40em height 0.5pt}}}; }}}\providecommand\AddOldversionToMargin[1]{\AddToHook {shipout/background}{\tikz[overlay,remember picture]{\node at ($(current page.west)+(2.5,0)$) [rotate=90] {\textcolor{\OldColor}{\vbox{\hrule width \the\textwidth height 0.5pt} \textcolor{black}{#1}\ \vbox{\hrule width 40em height 0.5pt}}}; }}}\providecommand\AddLineToMargin[3]{\AddToHook {shipout/background}{\tikz[overlay,remember picture]{\node at ($(current page.west)+(#2,0)$) [rotate=90] {\textcolor{#1}{\vbox{\hrule width \the\textwidth height 0.5pt} \textcolor{black}{#3}\ \vbox{\hrule width 40em height 0.5pt}}}; }}}\let\oldtocsection=\tocsection \let\oldtocsubsection=\tocsubsection \let\oldtocsubsubsection=\tocsubsubsection \renewcommand{\tocsection}[2]{\hspace{0em}\vspace*{0.1em}\oldtocsection{#1}{#2}}\renewcommand{\tocsubsection}[2]{\hspace{4ex}\oldtocsubsection{#1}{#2}}\renewcommand{\tocsubsubsection}[2]{\hspace{6ex}\oldtocsubsubsection{#1}{#2}}\usepackage{ulem}\usepackage{fancybox}\ifpdf \usepackage[pdftex]{lscape}\else \usepackage{lscape}\fi \makeatletter \newcommand{\verbatimfont}[1]{\def\verbatim@font{#1}}\makeatother \IfFileExists{mathabx.sty}{}{}\usepackage{amsfonts}\usepackage{amssymb}\usepackage{stmaryrd}\usepackage{amsmath}\usepackage{amsthm}\usepackage{dsfont}\IfFileExists{mbboard.sty}{\usepackage{mbboard}}{}\usepackage{mathrsfs}\usepackage{twcal}\usepackage{accents}\usepackage[T1]{fontenc}\usepackage[latin1]{inputenc}\catcode`\=13 \def{+}\newcommand\assigncharacter[1]{\expandafter\newcommand\csname #1\endcsname{\mathds{#1}}}\FunctionForEach{,}{\assigncharacter}{A,B,C,D,E,F,G,I,J,K,M,N,Q,R,T,U,V,W,X,Y,Z}\renewcommand\assigncharacter[1]{\expandafter\newcommand\csname C#1\endcsname{\mathcal{#1}}}\FunctionForEach{,}{\assigncharacter}{A,B,C,D,E,F,G,H,I,J,K,L,M,N,O,P,Q,R,S,T,U,V,W,X,Y,Z}\renewcommand\assigncharacter[1]{\expandafter\newcommand\csname D#1\endcsname{\mathfrak{#1}}}\FunctionForEach{,}{\assigncharacter}{a,b,c,d,e,f,g,h,i,j,k,l,m,n,o,p,q,r,s,t,u,v,w,x,y,z,A,B,C,D,E,F,G,I,K,L,M,N,O,P,Q,R,S,T,U,V,W,X,Y,Z} \renewcommand\assigncharacter[1]{\expandafter\newcommand\csname S#1\endcsname{\mathscr{#1}}}\FunctionForEach{,}{\assigncharacter}{A,B,C,D,E,F,G,H,I,J,K,L,M,N,O,P,Q,R,T,U,V,W,X,Y,Z}\def\NewFont#1#2#3#4#5{\expandafter\font\csname #1display\endcsname =#1 at #2 \expandafter\font\csname #1normal\endcsname =#1 at #3 \expandafter\font\csname #1script\endcsname =#1 at #4 \expandafter\font\csname #1scriptscript\endcsname =#1 at #5 }\def\NewFontLetter#1#2{{\mathchoice {{\expandafter\hbox{\csname #1display\endcsname\char"#2}}}{{\expandafter\hbox{\csname #1normal\endcsname\char"#2}}}{{\expandafter\hbox{\csname #1script\endcsname\char"#2}}}{{\expandafter\hbox{\csname #1scriptscript\endcsname\char"#2}}}}}\NewFont{pxsyc}{9.00pt}{8.00pt}{7.00pt}{6.00pt}\NewFont{pxsya}{9.00pt}{8.00pt}{7.00pt}{6.00pt}\NewFont{p1xr}{10.00pt}{9.00pt}{8.00pt}{7.00pt}\NewFont{MnSymbolA5}{10.00pt}{9.00pt}{8.00pt}{7.00pt}\NewFont{MnSymbolC10}{10.00pt}{9.00pt}{8.00pt}{7.00pt}\NewFont{MnSymbolD10}{12.00pt}{11.00pt}{10.00pt}{9.00pt}\NewFont{MnSymbolF10}{12.00pt}{11.00pt}{10.00pt}{9.00pt}\def\IndependenceX#1#2{#1\setbox0=\hbox{$#1x$}\kern\wd0\hbox to 0pt{\hss$#1\mid$\hss}\lower.9\ht0\hbox to 0pt{\hss$#1\smile$\hss}\kern\wd0}\def\nIndependenceX#1#2{#1\setbox0=\hbox{$#1x$}\kern\wd0 \hbox to 0pt{\mathchardef\nn=12854\hss$#1\nn$\kern1.4\wd0\hss}\hbox to 0pt{\hss$#1\mid$\hss}\lower.9\ht0 \hbox to 0pt{\hss$#1\smile$\hss}\kern\wd0}\NewFont{manfnt}{12.00pt}{11.00pt}{10.00pt}{9.00pt}\NewFont{favmr7y}{12.00pt}{11.00pt}{10.00pt}{9.00pt}\ifpdf \usepackage[pdftex,usenames,x11names]{xcolor}\else \usepackage[dvips,usenames,x11names]{xcolor}\fi \StackCreate{ColoR}{;}{?}\AtBeginDocument{\colorlet{defaultcolor}{.}\ColoRPush{defaultcolor}}\definecolor{Green}{rgb}{0.00,0.50,0.00}\definecolor{DarkGreen}{rgb}{0.00,0.40,0.00}\definecolor{gray}{rgb}{0.40,0.40,0.40} \renewcommand\textcolor[2]{\ColoRPush{#1}\color{\ColoRTop}#2\ColoRPop\color{\ColoRTop}}\usepackage[pdftex]{graphicx}\usepackage[all]{xy}\ifdefined\loadTIKZ \loadTIKZ \def\TIKZlabel#1{}\else \usepackage{tikz}\usetikzlibrary{matrix,arrows,calc,cd,decorations.pathmorphing,intersections}\fi \newcommand {\notion}[2][]{\def\temp{#1}\ifmmode #2 \ifx \temp\empty \index{$#2$}\else \index{$#1$}\fi \else {\bf #2}\ifx \temp\empty \index{#2}\else \index{#1}\fi \fi }\newcommand\NOPAGENUMBER[1]{}\usepackage{xr-hyper}\newcommand{\refX}[2]{\IfBeginWith{#1}{:}{\ref{\GetAfterColonAux#1;-#2}}{\cite[\ref{#1-#2}]{#1}}}\providecommand\hyperrefBOOKMARKS{true}\providecommand\hyperrefOPTIONS{destlabel=false}\ifdefined\BibLaTeX \usepackage[\hyperrefOPTIONS,pdftex,linktocpage,breaklinks,bookmarks=\hyperrefBOOKMARKS]{hyperref}\else \usepackage[\hyperrefOPTIONS,pdftex,linktocpage,pagebackref,breaklinks,bookmarks=\hyperrefBOOKMARKS]{hyperref}\fi \hypersetup{bookmarksdepth=3, colorlinks=true,allcolors=Green, linkcolor=DarkGreen, citecolor=violet, urlcolor=blue, runcolor=red, filecolor=magenta }\providecommand\url[1]{}\providecommand\nolinkurl[1]{}\providecommand\href[3][]{}\providecommand\hyperlink[2]{}\providecommand\hypertarget[2]{}\providecommand\hyperdef[3]{}\providecommand\hyperref[2][]{} \providecommand\hypersetup[1]{}\providecommand\pdfbookmark[3][]{}\providecommand\currentpdfbookmark[2]{}\providecommand\belowpdfbookmark[2]{}\providecommand\texorpdfstring[2]{} \def\UndefinedRef#1{\LARGE\bfseries\color{red} ??#1??}\makeatletter \def\@setref#1#2#3{\ifx#1\relax \protect\G@refundefinedtrue \nfss@text{\reset@font\UndefinedRef{#3}}\@latex@warning{Reference `#3' on page \thepage \space undefined }\else \expandafter\Hy@setref@link#1\@empty\@empty\@nil{#2}\fi }\makeatother \fi \newcommand\pr{\begin{proof}}\def\ende{\end{proof}}\newtheoremstyle{LayoutVoid}{1ex}{0ex}{\normalfont}{}{\bf}{.}{1ex}{}\newcommand\stressstatement[1]{#1}\theoremstyle{plain}\swapnumbers \newcommand\maketheorem[1]{\newtheorem{#1}[theorem]{\stressstatement{#1}} \newtheorem{#1Definition}[theorem]{\stressstatement{#1 and Definition}}  }\FunctionForEach{,}{\maketheorem}{Conclusion,Conjecture,Corollary,Fact,Facts,Lemma,Observation,Observations,Proposition,Reminder,Scholium,Summary,Theorem}\theoremstyle{definition}\theoremstyle{remark}\FunctionForEach{,}{\maketheorem}{Convention,Counterexample,Counterexamples,Discussion,Example,Examples, Exercise,Exercises,Explanation,Notation,Project,Projects,Question,Questions,Remark,Remarks,Strategy,Warning}\theoremstyle{LayoutVoid}\numberwithin{equation}{section}\newcommand{\labelon}[1]{\marginpar{#1}}\newcommand{\labelx}[1]{{\def\temp{#1}\ifx\temp\empty\else \label{#1}\labelon{#1}\fi}}\def\GetAfterColon#1:#2;;{#2}\def\GetAfterPlus#1+#2;;{#2}\newenvironment{FACT}[2]{\edef\currentlabel{#2}\IfBeginWith{#1}{:}{\def\tempFactName{void}\def\tempFreeTitle{\GetAfterColon#1;;\ }}{\IfBeginWith{#1}{+}{\def\tempFactName{voidTheorem}\def\tempFreeTitle{\GetAfterPlus#1;;\ }}{\def\tempFactName{#1}\def\tempFreeTitle{}}}\def\tempfacT{\end{\tempFactName}}\begin{\tempFactName}\labelx{#2}\textup{\textbf{\tempFreeTitle}}\capitalize[q]{#1}\caselower[q]{#1}}{\tempfacT}\newcommand{\inv}{\mathrm{inv}}\catcode95=12 \catcode95=8 \newcommand\bdl{{\ifmmode \mathrm{bdlat}\else {bounded distributive lattice}\fi}} \newcommand{\st}{{\ \vert \ }}\let\temp\phi \let\phi\varphi \let \varphi\temp \let\temp\theta \let\theta\vartheta \let \vartheta\temp \let\eps\varepsilon  \let\0\emptyset \newcommand{\into}{\hookrightarrow}\makeatletter \newcommand{\xRightarrow}[2][]{\ext@arrow 0359\Rightarrowfill@{#1}{#2}}\newcommand{\xLeftarrow}[2][]{\ext@arrow 0359\Leftarrowfill@{#1}{#2}}\newcommand{\xonto}[2][]{\ext@arrow 0359\rightarrowfill@ {#1}{#2}\mathrel{\mspace{-15mu}}\rightarrow}\newcommand{\xinto}[2][]{\lhook\joinrel\ext@arrow 0359\rightarrowfill@ {#1}{#2}}\makeatother \newcommand{\lra}{\longrightarrow}\newcommand{\ra}{\rightarrow}\newcommand{\Ra}{\Rightarrow}\newcommand{\mal}{{\cdot}} \newcommand{\addressTressl}{The University of Manchester, Department of Mathematics, Oxford Road, Manchester M13 9PL, UK}\newcommand{\emailTressl}{marcus.tressl@manchester.ac.uk}\newcommand{\homepageTressl}{\url{http://personalpages.manchester.ac.uk/staff/Marcus.Tressl/}}\providecommand{\cprime}{$'$}\newcommand{\monthname}[1]{\ifcase#1 \or January \or February \or March \or April \or May \or June \or July \or August \or September \or October \or November \or December \fi}\newcommand\LongColor{teal}\newcommand\OldColor{gray}\newcommand\COL{\ifmmode\colon\else :\ \fi}\newcommand{\claim}{\textit{Claim.}\ } \newcommand\operator[1]{\mathop{\operatorname{#1}}\nolimits} \newcommand{\Char}{\operator{char}}\newcommand{\GL}{\operator{GL}}\newcommand{\Gl}{\operator{GL}}\IfFileExists{C:/wb/System64/WinBatch.exe}{}{}\usepackage{todonotes}\renewcommand{\labelon}[1]{}\newcommand{\tr}{\operator{tr}}\renewcommand\inv{\rm invo}
\begin{document} \title{A Model Theoretic Perspective on Matrix Rings} \author{Igor Klep}\address{ Faculty of Mathematics and Physics, University of Ljubljana \& Faculty of Mathematics, Science and Information Technology, University of Primorska, Koper \& Institute of Mathematics, Physics and Mechanics, Ljubljana, Slovenia \newline Homepage: \url{https://igorklep.github.io/}} \email{igor.klep@fmf.uni-lj.si} \thanks{The first author is supported by the Slovenian Research Agency program P1-0222 and grants J1-50002, J1-2453, N1-0217, J1-3004. Part of this work prior to 2019 was partially supported by the Marsden Fund Council of the Royal Society of New Zealand} \author{Marcus Tressl}\address{\addressTressl\newline Homepage: \homepageTressl} \email{\emailTressl} \thanks{The second author was supported by MIMS, Department of Mathematics, The University of Manchester} \date{\today} \def\igor{\color{blue}} \begin{abstract} In this paper natural necessary and sufficient conditions for quantifier elimination of matrix rings $M_n(K)$ in the language of rings expanded by two unary functions, naming the trace and transposition, are identified. This is used together with invariant theory to prove quantifier elimination when $K$ is an intersection of real closed fields. On the other hand, it is shown that finding a natural \textit{definable} expansion with quantifier elimination of the theory of $M_n(\C)$ is closely related to the infamous simultaneous conjugacy problem in matrix theory. Finally, for various natural structures describing dimension-free matrices it is shown that no such elimination results can hold by establishing undecidability results. \end{abstract} \maketitle \setcounter{tocdepth}{3} \tableofcontents \section{Introduction} \normalem \par \noindent This article grew out of an attempt to understand the first order model theory of full matrix rings in connection with their use in what is called \textit{Free Analysis}. Free Analysis \cite{Voiculescu2010} studies functions in noncommuting variables (such as polynomials, rational functions or formal power series) and their evaluations in noncommutative algebras, such as matrix rings or operators on Hilbert space, and thus provides a framework for dealing with quantities with the highest degree of non-commutativity, such as large random matrices, see for example \cite{AglMcC2016,KaVVin2014,HeKlMc2011}. Our focus lies in the \textit{free algebra} consisting of polynomials in noncommuting variables (=noncommutative polynomials) over $\K $ ($=\R$ or $\C$). Here, noncommutative polynomials are considered as functions by evaluating them in matrix rings (of arbitrary size!). For example the polynomial $xy-yx$ can only be distinguished from $0$ by evaluating at matrices, say of size $2\times 2$. By the celebrated Amitsur-Levitzki theorem \cite{rowen,procesiBook}, if $P({\bar x})$ is an arbitrary noncommutative polynomial then $P$ is 0 if and only if $P({\bar X})=0$ for all tuples ${\bar X}$ of square matrices \textit{of any size}. Thus, matrices in this context are not restricted to a specific size and we may refer to ``dimension-free'' matrices when we want to stress this point of view. \par There is ongoing interest (cf. \cite{DrNeTh2023,Putinar2007}) in the question of whether some form of elimination theory or decidability from the classical case of the field $\K$ can be rescued in the noncommutative context. This article contributes to these questions in two ways. To explain how, first note that these questions have obvious negative answers if we ask them for the common theory of all matrix rings $M_n(\K)$; this theory is not model-complete in the language of rings as $M_n(\K)$ is not elementary in $M_{n+1}(\K)$ for any $n$, and it is indeed hereditarily undecidable because every non-principal ultraproduct of the $M_n(\K)$ interprets true arithmetic. \par A less na\"ive way to tackle the problem is to consider first order structures that interpret all $M_n(\K)$ and then to try to approach elimination theory and decidability questions for such structures. We show that the most commonly used structures in Free Analysis, which interpret all matrix rings $M_n(\K)$, are undecidable. This is done in Section \ref{SectionUndecidable}. (We point out that the community has not agreed on the exact structure to be used for dimension-free matrices yet.) This already implies that quantifier elimination results similar to those for algebraically closed fields or real closed fields cannot be expected to hold true for structures interpreting all $M_n(\K)$. However it is unclear if a weakened elimination result like model-completeness holds true in a suitable language. \par In this context it is important to understand the elimination theory of matrix rings of \textit{fixed size}, which, surprisingly, is strongly tied to the well-known simultaneous conjugacy problem for matrices (asking for invariants, or, normal forms, for the conjugacy class of a pair of $n\times n$-matrices under the action of $\GL_n$; this problem lies at the bottom of the tame-wild dichotomy in the representation theory of finite dimensional algebras, see \cite[p.~vii, last paragraph]{BaSaZu2009}). To be more precise, let $\SL$ be the first-order language of rings. The question about the elimination theory of $M_n(\K)$ in this language, a priori, seems to be all answered by the classical results for the field $\K$. (Contemporary model theory might even identify the bi-interpretable structures $\K$ and $M_n(\K)$.) However, already $M_n(\C)$ does not have quantifier elimination in $\SL$ (cf. \ref{NoParamNoQE}) and it admits quantifier elimination only if invariants for the simultaneous conjugacy problem are named in an extended language. This is done in Section \ref{QEoverC}.\looseness=-1 \par We now explain our main contribution, namely the elimination theory of matrix rings of fixed size $n\times n$. We switch to an arbitrary field $K$. A classical comparison of the field $K$ and the matrix ring $M_n(K)$ in terms of \textit{how} the bi-interpretation is done reveals a more subtle elimination theory of $M_n(K)$. For an example, consider polynomials $P(x,y),Q(x,y)$. The solution set in $M_n(K)^2$ of $P(x,y)=0, Q(x,y)\neq 0$, seen as a subset of $K^{2\mal n^2}$, is closed under simultaneous conjugation. The question of whether the projection onto the $X$-coordinate(s) has this property is not answered within the elimination theory of $K$. The issue is that the quantifier-free definable sets in $M_n(K)$ (in the language of rings for now) single out certain $K$-definable sets and not all $K$-varieties can be described quantifier-free in $M_n(K)$. The ring $M_n(K)$ is quantifier-free definable in the field $K$. Conversely, $K$ is universally definable in the ring $M_n(K)$ as its center\footnote{It should also be noted that for any field $K$, the ring $M_n(K)$ is already interpretable in the monoid $(M_n(K),\cdot)$ when $n\geq 3$. The reason is that $(M_n(K),\cdot)$ interprets the poset of vector subspaces of $K^n$ and one can then invoke incidence geometry, see \cite[5.1]{Tressl2017}. For the interpretation we code a subspace as the range of a matrix and note that $\mathrm{ran}(A)\subseteq \mathrm{ran}(B)\iff \exists C\in M_n(K): A=BC$.} and in \ref{CenterExDef} we see a positive primitive (in particular, existential) definition. However there is no field $K$ that is quantifier-free definable in the ring $M_n(K)$ as its center, see \ref{SectionExamples}. \par In Section \ref{SectionQE} we identify natural necessary and sufficient conditions for quantifier elimination of $M_n(K)$ in the language of rings expanded by two unary functions, naming the trace and transposition. This is obtained for formally real Pythagorean fields (see \ref{MnTrHqe}) and it says that $M_n(K)$ has quantifier elimination in the extended language if and only if there is some $D\in \N$ depending only on $n$ such that for all $d$ and any two $d$-tuples of $n\times n$ matrices $X,Y\in M_n(K)^d$ with \[ \tr(w(X,X^t))=\tr(w(Y,Y^t)) \] for all words $w$ in $x,x^t$ of length $\leq D$, there is some $O\in M_n(K)$ with $OO^t=I_n$ and $O^tX_iO=Y_i$ for all $i$, i.e., the tuples $X$ and $Y$ are orthogonally equivalent over $K$.\looseness=-1 \par This condition is satisfied for the field of real numbers and more generally for every intersection of real closed fields, see \ref{prop:Specht}. A similar result holds for the complex field, however the involution properly expands the matrix ring to include the reals. As mentioned above, quantifier elimination of a natural \textit{definable} expansion of $M_n(\C)$ is closely related to the simultaneous conjugacy problem; see \ref{QEoverC}. \par \medskip For the theory of matrix rings and more generally, C$^*$-algebras from a continuous logic perspective we refer the reader to e.g. \cite{FaHaSh2014} (notice however, that our goals and our results are tightly bound to the first order theory of matrix rings). \looseness=-1 We use basic model theory and standard notations as explained for example in \cite{Hodges1993}. For generalities on decidability in first order logic see \cite{Rauten2010}. All rings and algebras in this paper are associative but not necessarily commutative or unital, i.e., they might not have an identity element. Fields are commutative. \section{Elimination theory of matrix rings} \labelx{SectionQE} \numberwithin{theorem}{subsection} \noindent In this section we are concerned with the elimination theory of matrix rings of fixed size. The first subsection is of preliminary nature and deals with model-completeness in the ring language and with quantifier elimination after naming matrix units. After that, in the main part, we study natural expansions by trace and transposition (or adjoint). In particular, we prove quantifier elimination of the ring $M_n(\R)$ expanded by the trace, transposition and the order on its center, see \ref{cor:PAP}.\looseness=-1 \subsection{Naming matrix units} In this subsection we show that model-complete expansions of fields have model-complete matrix rings in their natural language, see \ref{MnKmodelcomplete}. If we name matrix units, the same is true for quantifier elimination, see \ref{MnKMatUnitsQE}.\looseness=-1 \begin{FACT}{:On matrix units.}{MatrixEijElementary} We describe the abstract properties of a full set of matrix units (i.e., $n\times n$-matrices that have exactly one entry 1 and all other entries 0). Let $A$ be a ring, $n\in \N$ and for $i,j\in \{1,\ldots,n\}$ let $a_{ij}\in A$. Suppose for all $i,j,s,t\in \{1,\ldots,n\}$ we have $ a_{ij}\mal a_{st}= \delta_{js} a_{it} $. \noindent The following properties are easily verified.\looseness=-1 \begin{enumerate}[\rm(1),itemsep=1ex] \item For $i,j,s,t\in \{1,\ldots,n\}$ we have $a_{ss}a_{ij}a_{tt}=\delta_{is}\delta_{jt}a_{ij}$. \item If $a_{ij}=0$ for some $i,j$, then $a_{st}=a_{si}\mal a_{ij}\mal a_{jt}=0$ for all $s,t$. Now assume all $a_{ij}\neq 0$. Then the $a_{ij}$ ($1\leq i,j\leq n$) are linearly independent over any central subfield $F$ of $A$. \item For $(x_{ij})_{i,j\in\{1,\ldots,n\}},\ (y_{ij})_{i,j\in\{1,\ldots,n\}}\in M_n(F)$, we have \[ (\sum_{i,j=1}^nx_{ij}a_{ij})\mal (\sum_{i,j=1}^ny_{ij}a_{ij})= \sum_{i,j=1}^n(\sum_{k=1}^nx_{ik}y_{kj})a_{ij}. \] \item Let $F$ be a central subfield of $A$. The map \begin{align*} M_n(F)\lra A, \quad (u_{ij})_{i,j\in\{1,\ldots,n\}}\longmapsto \sum_{i,j=1}^nu_{ij}a_{ij} \end{align*} is a (not necessarily unital) $F$-algebra homomorphism, because it is clearly $F$-linear and it is a ring homomorphism by (3). If $a_{ij}\neq 0$ for all $i,j\in\{1,\ldots,n\}$, then by (2) this map is injective. \par To see an example where the map is not unital, choose any field $F$, set $n=1<m$, $A=M_m(F)$ and take $a_{11}\in A\setminus \{0,I_m\}$ with $a_{11}^2=a_{11}$. \end{enumerate} \end{FACT} \begin{FACT}{:Defining matrix units.}{} The language of unital rings is denoted by $ \SL_\mathrm{ri}=\{+,\cdot,-,0,1\}. $ Let $F$ be a field and let $M=M_n(F)$. The center $C=C_n$ of $M$ is isomorphic to $F$, but we will work with $C$ instead of $F$. For $N\in \N$ we consider $M_N(C)$ as a subset of $M^{N^2}$ and as an $F$-algebra via the natural embedding $F\cong C\into M_N(C)$. Take $2N^2+2$ variables $ {\bar u}=(u_{ij}\st i,j\in \{1,\ldots N\}),\ {\bar x}=(x_{ij}\st i,j\in \{1,\ldots N\}),\ y,\ v. $ Consider the following $\SL_\mathrm{ri}$-formulas: \begin{enumerate}[\rm(1),leftmargin=4ex] \item Let $\eps=\eps_N({\bar u})$ be the formula \[ \bigwedge_{i,j,t=1}^Nu_{ij}\mal u_{jt}=u_{it}\neq 0\ \land\ \bigwedge_{i,j,s,t=1, j\neq s}^Nu_{ij}\mal u_{st}=0 \] {\color{black}and let $\eps^+=\eps^+_N({\bar u})$ be the formula \[ \bigwedge_{i,j,t=1}^Nu_{ij}\mal u_{jt}=u_{it}\ \land\ \sum_{i=1}^Nu_{ii}=1\ \land\bigwedge_{i,j,s,t=1, j\neq s}^Nu_{ij}\mal u_{st}=0. \]} \item Let $\delta=\delta_N(v,{\bar u})$ be the formula $ \bigwedge_{s,t=1}^Nv\mal u_{st}=u_{st}\mal v. $ \item Let $\lambda_N({\bar x},y,{\bar u})$ be the formula $ y=\sum_{i,j=1}^Nx_{ij}\mal u_{ij}. $ \end{enumerate} Finally let $\gamma=\gamma_N({\bar x},y,{\bar u})$ be the formula $ \lambda({\bar x},y,{\bar u})\land \eps({\bar u})\land \bigwedge_{i,j=1}^N\delta(x_{ij},{\bar u}). $ \end{FACT} \noindent By \ref{MatrixEijElementary} we then obtain \begin{FACT}{Proposition}{MatrixEliminate} For $i,j\in \{1,\ldots,N\}$ let $E_{ij}\in M_N(C)$ be the $N\times N$-matrix that has exactly one nonzero entry, namely $1$ ($\in C$) at position $(i,j)$. \begin{enumerate}[\rm(1)] \item If $\Theta :M_N(C)\lra M_n(F)=M$ is a (not necessarily unital) embedding of $F$-algebras, then the $N^2$-tuple ${\bar a}:=(\Theta (E_{ij}))_{i,j\in \{1,\ldots N\}}\in M^{N^2}$ is a realization of $\eps_N({\bar u})$, and $\gamma_N({\bar x},y,{\bar a})$ defines the graph of $\Theta $ in the ring $M$. \item For every realization ${\bar a}=(a_{ij})_{i,j\in \{1,\ldots,N\}}\in M^{N^2}$ of $\eps_N$ in $M$, there is a unique (not necessarily unital) embedding of $F$-algebras $\Theta_{\bar a} :M_N(C)\lra M_n(F)$ such that $\Theta_{\bar a} (E_{ij})=a_{ij}$ $(i,j\in \{1,\ldots,N\})$. Explicitly, the graph of $\Theta_{\bar a}$ is defined by $\gamma_N({\bar x},y,{\bar a})$. \end{enumerate} Consequently, \begin{enumerate}[resume*] \item The family of all (not necessarily unital) embeddings of $F$-algebras $M_N(C)\lra M_n(F)$ is quantifier-free definable in $M$ by $\gamma({\bar x},y,{\bar u})$ and its parameter set is quantifier-free defined by $\eps({\bar u})$. \item{\color{black}The formulas $\eps_N$ and $\eps_N^+$ define the same set in the ring $M_N(C)$.}\qed \end{enumerate} \end{FACT} \begin{FACT}{Corollary}{CenterExDef} For any field $F$ the center of $M_n(F)$ is positive primitively definable in the language $\SL_\mathrm{ri}$ by \( \exists {\bar u}\bigl(\eps_n^+({\bar u})\land\ \delta_n(v,{\bar u})\bigr). \) \qed \end{FACT} \begin{FACT}{Corollary}{MnAxiom} \begin{enumerate}[\rm(1)] \item\label{it:mn1} For a field $F$, the theory of $M_n(F)$ is axiomatised by saying the following about a model $A$ with center $C$: \begin{enumerate}[\rm (a)] \item $A$ is a ring whose center $C$ is elementarily equivalent to $F$. \item There is some realization ${\bar a}=(a_{ij})_{i,j\in \{1,\ldots,n\}}$ of $\eps^+_n$ in $A^{n^2}$ and for each such realization, $\gamma_n({\bar x},y,{\bar a})$ defines an isomorphism $M_n(C)\lra A$. \end{enumerate} \item\label{it:mn2} If $A,B$ are rings that are elementarily equivalent to $M_n(F)$, and if $A$ is a subring of $B$, then the center $C_A$ of $A$ is a subring of $C_B$. Further, for each realization ${\bar a}=(a_{ij})_{i,j\in \{1,\ldots,n\}}$ of $\eps^+_n$ in $A^{n^2}$ the following diagram commutes: \begin{center} \begin{tikzcd}[row sep=8ex,column sep=8ex] A \ar[r,hook] & B \\ M_n(C_A) \ar[u, "\Theta_{\bar a}", "\cong"'] \ar[r,hook] & M_n(C_B) \ar[u,"\Theta_{\bar a}"',"\cong"] \end{tikzcd} \end{center} \vspace{-5mm}\qed \end{enumerate} \end{FACT} \begin{FACT}{Definition}{} Let $F$ be a field and let $\tilde F$ be an expansion of $F$ in some language $\SL$ extending $\SL_\mathrm{ri}$ (cf. \cite[p.9]{Hodges1993}). Then we define the $\SL$-structure $M_n(\tilde F)$ as the structure expanding the ring $M_n(F)$ and that interprets new relation symbols and constant symbols only on the center $C$ of $M_n(F)$ as given by $\tilde F$. A new $m$-ary function symbol $f$ is interpreted on $C^m$ as given by $\tilde F$, and set to be $0$ outside of $C^m$.\looseness=-1 \end{FACT} \noindent Recall that a structure $M$ in some language $\SL$ is called model-complete if the $\SL$-theory of $M$ is model-complete (cf. \cite[Thm 8.3.1]{Hodges1993}). Similarly, $M$ has quantifier elimination, if its $\SL$-theory has quantifier elimination. \par \smallskip\noindent For algebraically closed fields, the following may be found in \cite[Theorem 5.4]{Rose1980a}. \begin{FACT}{Proposition}{MnKmodelcomplete} If $\tilde F$ is a model-complete expansion of a field $F$ in some language $\SL$ extending $\SL_\mathrm{ri}$, then the $\SL$-structure $M_n(\tilde F)$ is also model-complete. Hence, for example, the ring $M_n(\C)$ is model-complete and the ring $M_n(\R)$ expanded by the natural order on its center is model-complete. \end{FACT} \begin{proof} This is a routine argument using \ref{MnAxiom}: Let $\tilde A,\tilde B$ be $\SL$-structures with underlying rings $A,B$ respectively. Suppose $\tilde A,\tilde B$ are elementarily equivalent to $M_n(\tilde F)$ with $\tilde A\subseteq \tilde B$. We need to show that $\tilde A\prec \tilde B$. Choose a realization ${\bar a}=(a_{ij})_{i,j\in \{1,\ldots,n\}}$ of $\eps^+_n$ in $A^{n^2}$ as in \ref{MnAxiom}\ref{it:mn1} and consider the commutative diagram of \ref{MnAxiom}\ref{it:mn2}. We see that the $\SL$-structure $\SM$ induced by $\tilde A$ on $C_A$ is a substructure of the $\SL$-structure $\SN$ induced by $\tilde B$ on $C_B$. By assumption this extension is elementary. Since $\tilde A$ is interpretable in $\SM$ in the same way $\tilde B$ is interpretable in $\SN$, we get $\tilde A\prec\tilde B$. \end{proof} \begin{FACT}{Remark}{NoParamNoQE} A corresponding version of \ref{MnKmodelcomplete} for quantifier elimination (instead of model-completeness) fails; for instance the ring $M_n(\C)$ does not have quantifier elimination in $\SL_\mathrm{ri}$ for any $n\geq 2$. In fact, by \cite[proof of Theorem 3.2]{Rose1978a}, for any infinite field $F$, the center of $M_n(F)$ is not quantifier-free definable with parameters from $F\mal I_n$ in the ring $M_n(F)$. \par A geometric argument goes as follows: Assume $F\cdot I_n$ is quantifier-free $F\cdot I_n$-definable in $M_n(F)$. Then $F\cdot I_n$ is a finite union of nonempty sets of the form $\{X\in M_n(F) \mid p_1(X)=\ldots=p_r(X)=0 \text{ and } q_1(X),\ldots,q_s(X)\neq 0\},$ where $p_i,q_j$ are univariate polynomials from $F[t]$. Since such polynomials have only finitely many roots in $F$ and $F$ is infinite, one of these sets is of the form $\{X\in M_n(F) \mid q_1(X),\ldots,q_s(X)\neq 0\}.$ But then $F\mal I_n$ has nonempty Zariski interior in $M_n(F)$, a contradiction. \end{FACT} \noindent If we allow matrix units as parameters, then a corresponding version of \ref{MnKmodelcomplete} for quantifier elimination does hold. \begin{FACT}{Lemma}{MatUnitsPresent} If $U$ is a subring of $M_n(F)$, $F$ a field and $U$ contains the standard matrix units $E_{ij}$, $1\leq i,j\leq n$, then \[ R_U=\{a\in F\st a\text{ is the $(1,1)$ entry of some }Y\in U\} \] is a subring of $F$ and $U=M_n(R_U)$. \end{FACT} \begin{proof} Let $a,b\in R_U$, say $a$ is the $(1,1)$ entry of $X\in U$, and $b$ is the $(1,1)$ entry of $Y\in U$. Then $a+b$ is the $(1,1)$ entry of $X+Y\in U$, and $ab$ is the $(1,1)$ entry of $XE_{11}YE_{11}\in U$, proving then $R_U$ is a subring of $k$. \par Given $X\in U\subseteq M_n(F)$, $X=(x_{ij})_{i,j}$, we see that $x_{ij}$ is the (1,1) entry of $E_{1i}XE_{j1}\in U$, so $U\subseteq M_n(R_U)$. Conversely, if $X=(x_{ij})_{i,j}\in M_n(R_U),$ then each $x_{ij}$ is the (1,1) entry of some $Y_{ij}\in U$. Hence \( X=\sum_{i,j} E_{i1} Y_{ij} E_{1j} \in U. \) \end{proof} \begin{FACT}{void}{QEisMCplusAP} Recall from \cite[Thm. 8.4.1]{Hodges1993} that an $\SL$-theory $T$ has quantifier elimination if and only if it is model-complete and models of $T$ have the amalgamation property over substructures. \end{FACT} \begin{FACT}{Proposition}{MnKMatUnitsQE} Let $\tilde F$ be an expansion with quantifier elimination of a field $F$ in some language $\SL$ extending $\SL_\mathrm{ri}$ and let ${\bar c}=(c_{ij})_{i,j\in \{1,\ldots,n\}}$ be new constant symbols. Then the $\SL({\bar c})$-structure $(M_n(\tilde F),{\bar e})$, where ${\bar c}$ is interpreted by a tuple ${\bar e}$ of matrix units, also has quantifier elimination. \par In particular, the ring $M_n(\C)$ expanded by the standard matrix units $E_{ij}$ and the ring $M_n(\R)$ expanded by the natural order on its center and the standard matrix units $E_{ij}$ have quantifier elimination. \end{FACT} \begin{proof} Since $M_n(\tilde F)$ is model-complete by \ref{MnKmodelcomplete}, it suffices to show that the theory of $(M_n(\tilde F),{\bar e})$ has the amalgamation property. Let $(\tilde A,{\bar a}),(\tilde B,{\bar b})$ be $\SL({\bar c})$-structures with underlying rings $A,B$ respectively. Suppose $(\tilde A,{\bar a}),(\tilde B,{\bar b})$ are elementarily equivalent to $(M_n(\tilde F),{\bar e})$ and suppose $\SU$ is a common $\SL({\bar c})$-substructure. Hence $\SU=(\tilde U,{\bar u})$, where $\tilde U$ is an expansion of a common subring $U$ of $A$ and $B$, and ${\bar u}={\bar a}={\bar b}$. Let $K,L$ be the center of $A,B$ respectively. By \ref{MnAxiom} there are ring isomorphisms $\phi:A\lra M_n(K),\psi:B\lra M_n(L)$ that map $u_{ij}$ to the standard matrix unit $E_{ij}$ for all $i,j$. We expand $M_n(K)$ to the $\SL$-structure $M_n(\tilde K)$ that makes $\phi$ an $\SL$-isomorphism $\tilde A\lra M_n(\tilde K)$, and similarly for $M_n(L)$. By \ref{MatUnitsPresent}, there are subrings $R\subseteq K,S\subseteq L$ such that the restriction of $\phi,\psi$ to $U$ are isomorphisms onto $M_n(R),M_n(S)$ respectively. We expand $M_n(R),M_n(S)$ to the induced $\SL$-substructures of $M_n(\tilde K),M_n(\tilde L)$ respectively and obtain the following commutative diagram:\looseness=-1 \begin{center} \begin{tikzcd} (M_n(\tilde K),{\bar E}) & (\tilde A,{\bar u}) \ar[l,"\phi","\cong"'] & \ & (\tilde B,{\bar u}) \ar[r,"\psi"',"\cong"] & (M_n(\tilde L),{\bar E})\\ (M_n(\tilde R),{\bar E})\ar[u,hook] & \ & (\tilde U,{\bar u}) \ar[ul, hook]\ar[ur, hook']\ar[ll,"\cong"']\ar[rr,"\cong"] & \ & (M_n(\tilde R),{\bar E})\ar[u,hook] \end{tikzcd} \end{center} Restricting all maps to centers and using that $\tilde F$ has quantifier elimination, there is some $\tilde \Omega$ elementarily equivalent to $\tilde K$ and $\tilde L$ together with $\SL$-embeddings $\rho:\tilde K\lra \tilde \Omega$, $\delta:\tilde L\lra \tilde \Omega$ such that for every $v$ in the center of $U$ we have $\rho(\phi(v))=\delta(\psi(v))$. Let ${\bar \rho}:M_n(K)\lra M_n(\Omega)$, ${\bar \delta}:M_n(L)\lra M_n(\Omega)$ be the unique extensions of $\rho,\delta$ preserving the standard matrix units. We see that ${\bar \rho},{\bar \delta}$ are $\SL({\bar c})$-morphisms and thus the desired amalgamation is given by the maps ${\bar \rho}\circ \phi$ and ${\bar \delta}\circ \psi$. \end{proof} \subsection{Quantifier elimination with trace and transposition} We have seen in \ref{MnKMatUnitsQE} that quantifier elimination of a field in a suitable language carries over to its matrix rings if we allow naming of definable parameters (i.e., the set of these parameters is 0-definable). Without parameters the assertion fails, see \ref{NoParamNoQE}. We now consider quantifier elimination of expansions of matrix rings by trace and transposition in the case of Pythagorean fields. We will see in \ref{MnTrHqe} that quantifier elimination is equivalent to a property in invariant theory describing simultaneous orthogonal similarity of matrices (where the conjugating matrix is orthogonal). For the real field the characterization entails quantifier elimination of the ring $M_n(\R)$ expanded by the trace, transposition and the order on its center. \begin{FACT}{Lemma}{AmalgTraces} Let $K,L$ be fields. Let $\SL$ be the extension of $\SL_\mathrm{ri}$ by a unary function symbol. Consider the $\SL$-structures $(M_n(K),\tr_K)$ and $(M_n(L),\tr_L)$. Let $(U,f)$ be an $\SL$-structure and suppose we are given $\SL$-embeddings $\phi:(U,f)\into (M_n(K),\tr_K)$ and $\psi:(U,f)\into (M_n(L),\tr_L)$. Then \begin{enumerate}[\rm(1)] \item\label{it:amal1} The subring $R$ of $U$ generated by the image of $f$ is commutative and $\phi(R)\subseteq K\mal I_n$, $\psi(R)\subseteq L\mal I_n$. \item\label{it:amal2} If $K\mal I_n$ and $L\mal I_n$ can be amalgamated over $\phi|_R,\psi|_R$ into some field $\Omega$ by maps $\rho:K\mal I_n\lra \Omega\mal I_n,\delta:L\mal I_n\lra \Omega\mal I_n$, then for the induced maps ${\bar \rho}:M_n(K)\lra M_n(\Omega), {\bar \delta}:M_n(L)\lra M_n(\Omega)$ and every $X\in U$ we have \[ \tr_\Omega({\bar \rho}(\phi(X)))=\tr_\Omega({\bar \delta}(\psi(X))). \] \noindent Here are the maps in a (not necessarily commutative) diagram. 
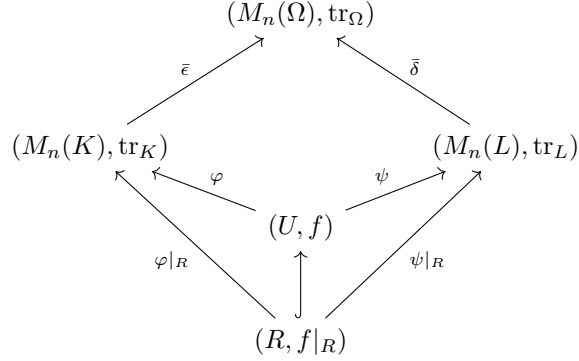
\begin{center} \begin{tikzcd}[row sep=3ex,column sep=3ex] & (M_n(\Omega),\tr_\Omega) \\ {} \\ (M_n(K),\tr_K) \arrow[uur, "\bar \rho"] && (M_n(L),\tr_L) \arrow[uul, "\bar \delta"'] \\ & (U,f) \arrow[ul,"\phi"'] \arrow[ur,"\psi"] \\ \\ & (R,f|_R) \arrow[uuul,"\phi|_R"] \arrow[uuur,"\psi|_R"'] \arrow[uu, hook] \end{tikzcd} \end{center} \end{enumerate} \end{FACT} \begin{proof} (1) Let $X\in U$, then $\phi(f(X))=\tr_K(\phi(X))$ since $\phi$ is an $\SL$-homomorphism $(U,f)\lra (M_n(K),\tr_K)$. Since $\tr_K(\phi(X))\in K\mal I_n$ we get $\phi(f(X))\in K\mal I_n$. Hence $\phi(f(U))\subseteq K\mal I_n$. Since $\phi$ is an embedding $U\lra M_n(K)$, $R$ is commutative and $\phi(R)\subseteq K\mal I_n$. Similarly, $\psi(R)\subseteq L\mal I_n$. \par \smallskip\noindent (2) For $X\in U$ we have \begin{align*} \tr_\Omega({\bar \rho}(\phi(X)))&=\rho (\tr_K(\phi(X)))\text{ since }\tr_\Omega\circ \,{\bar \rho}=\rho \circ \tr_K\cr &=\rho(\phi(f(X)))\text{ since }\tr_K\circ\,\phi= \phi\circ f\cr &=\delta(\psi(f(X)))\text{ since }\rho\circ \phi=\delta\circ \psi, \end{align*} and similarly $\tr_\Omega({\bar \delta}(\psi(X)))=\delta(\psi(f(X)))$. \end{proof} \begin{FACT}{Theorem}{Specht} Let $\Omega$ be a real closed field or the algebraic closure of a real closed field. For $X_1,\ldots,X_d,Y_1,\ldots,Y_d\in M_n(\Omega)$ the following are equivalent: \begin{enumerate}[\rm(1)] \item There is some unitary $O\in M_n(\Omega)$ with $O\mal X_i\mal O^*=Y_i$ for all $i\in\{1,\ldots,d\}$.\footnote{If $\Omega$ is real closed then $X^*$ is the transpose of $X$. If $\Omega$ is the algebraic closure of a real closed field $\Omega_0\subseteq \Omega$ then $X^*$ is the conjugate transpose of $X$ with respect to $\Omega_0$.} \item\label{it:allWords} For every word $w$ in the letters $x_1,\ldots,x_d,x_1^*,\ldots,x_d^*$ we have \[ \tr_\Omega(w(X_1,\ldots,X_d,X_1^*,\ldots,X_d^*))=\tr_\Omega(w(Y_1,\ldots,Y_d,Y_1^*,\ldots,Y_d^*)). \] \item For every word $w$ of degree $\leq n^2$ in the letters $x_1,\ldots,x_d,x_1^*,\ldots,x_d^*$, \[ \tr_\Omega(w(X_1,\ldots,X_d,X_1^*,\ldots,X_d^*))=\tr_\Omega(w(Y_1,\ldots,Y_d,Y_1^*,\ldots,Y_d^*)). \] \end{enumerate} \end{FACT}\begin{proof} The equivalence of (1) and (2) over $\C$ is established in \cite[Thm. 4]{Wiegma1961} and in \cite[Cor. 1]{Sibir1968}. The equivalence of (1) and (2) over $\R$ is given by \cite[Lemma 2]{Sibir1968} (see also \cite[Thm 7.1, Thm. 15.3]{Procesi1976}). For degree bounds in (3) (when $\Omega=\R$ or $\C$), see \cite[Thm 7.3]{Procesi1976} and \cite{Razmys1974}.\footnote{For $d=1$, this result is classical. The equivalence between (1) and (2) over $\C$ is due to \cite[Satz 1]{Specht1940}. The degree bounds and the real case for $d=1$ are due to \cite[Thm. 1 and Cor. to Thm. 2]{Pearcy1962}.} Since (2) is equivalent to (3), all equivalences carry over to all real closed fields and to their algebraic closures. \end{proof} \begin{FACT}{Observation}{lem:CCt} Let $F$ be a formally real field. Then \[ X=0 \iff \tr(X^tX)=0 \] for every matrix $X=(x_{ij})\in M_n(F)$, because $\tr(X^tX)=\sum_{i,j}x_{ij}^2$. \end{FACT} \begin{FACT}{Theorem}{MnTrHqe} Let $F$ be a formally real Pythagorean field (hence sums of squares are squares) and let $\tilde F$ be an expansion of $F$ in a language $\SL$ extending the language $\SL_\mathrm{ri}$. Suppose $\tilde F$ has quantifier elimination in $\SL$. Let $\SL(\tr,\inv)$ be the extension of $\SL$ by two new unary function symbols. The following are equivalent. \begin{enumerate}[\rm (1)] \item The structure $(M_n(\tilde F),\tr_F, X\mapsto X^t)$ has quantifier elimination in $\SL(\tr,\inv)$. \item $F$ has the Specht property for the transpose, i.e., there is some $D=D(n)$ such that for all $d$ and any two $d$-tuples of $n\times n$ matrices $X,Y\in M_n(F)^d$ with \[ \tr(w(X,X^t))=\tr(w(Y,Y^t)) \] for all words $w$ in $x,x^t$ of length $\leq D$, there is some $O\in M_n(F)$ with $OO^t=I_n$ and $O^tX_iO=Y_i$ for all $i$. \item If $\tilde K\equiv \tilde F$ and $\SU$ is a substructure of $(M_n(\tilde K),\tr_K,X\mapsto X^t)$ and $\psi:\SU\lra (M_n(\tilde K),\tr_K,X\mapsto X^t)$ is an embedding, then there is an elementary extension $\tilde \Omega \succ \tilde K$ and an extension of $\psi$ to an embedding $(M_n(\tilde K),\tr_K,X\mapsto X^t)\lra (M_n(\tilde \Omega),\tr_\Omega,X\mapsto X^t)$. Hence the following diagram commutes: \begin{center} \begin{tikzcd}[row sep=5ex,column sep=1ex] \ & (M_n(\tilde \Omega),\tr_\Omega,X\mapsto X^t) & \ \\ (M_n(\tilde K),\tr_K,X\mapsto X^t) \ar[ur] & \ & (M_n(\tilde K),\tr_K,X\mapsto X^t)\ar[ul,"\succ "', hook] \\ \ & \SU \ar[lu,hook]\ar[ru,"\psi"] & \ \end{tikzcd} \end{center} \end{enumerate} \end{FACT} \begin{proof} (2)$\Ra$(1) Since $\tilde F$ is model-complete we know from \ref{MnKmodelcomplete} that $M_n(\tilde F)$ is model-complete and so is its definable expansion $(M_n(\tilde F),\tr_F, X\mapsto X^t)$. Hence by \ref{QEisMCplusAP} it suffices to show that the theory $T$ of $(M_n(\tilde F),\tr_F, X\mapsto X^t)$ has the amalgamation property over finitely generated substructures. So let $\SM,\SN\models T$ and let $\SU$ be a common finitely generated $\SL(\tr,\inv)$-substructure of $\SM,\SN$. Using \ref{MatrixEijElementary}, \ref{MatrixEliminate} and as $\SM\models T$ we see that there is an isomorphism ${\bar \phi }:\SM\lra (M_n(\tilde K),\tr_K,X\mapsto X^t)$ where $\tilde K\equiv \tilde F$: In the language $\SL(\tr,\inv)$ we can say that there are matrix units $a_{ij}$ over the center $K$ of $\SM$ such that the ring homomorphism $M_n(K)\lra \SM$ that maps $E_{ij}$ to $a_{ij}$, is an isomorphism mapping transposition to the action of $\inv^\SM$.\looseness=-1 \par We write $\phi$ for the restriction of ${\bar \phi }$ to $\SU$. Similarly, we see that there is an isomorphism ${\bar \psi }:\SN\lra (M_n(\tilde L),\tr_L,X\mapsto X^t)$, with $\tilde L\equiv \tilde F$ and we write $\psi$ for the restriction of ${\bar \psi }$ to $\SU$. We now replace $\SM$ by $(M_n(\tilde K),\tr_K,X\mapsto X^t)$ and $\SN$ by $(M_n(\tilde L),\tr_L,X\mapsto X^t)$ and we need to amalgamate these $\SL(\tr,\inv)$ structures over $\SU$ via the $\SL(\tr,\inv)$-embeddings $\phi,\psi$. We write $\SU=(\tilde U,f,h)$, where $f:U\lra U$ and $h:U\lra U$ are the maps induced by the trace functions and the transpositions, respectively on $U$. \par Let $R$ be the subring of $U$ generated by the image of $f$. By \ref{AmalgTraces}\ref{it:amal1}, $R$ is commutative, $\phi(R)\subseteq K\mal I_n$ and $\psi(R)\subseteq L\mal I_n$. For better readability we now identify $K$ with $K\mal I_n$ and $L$ with $L\mal I_n$. Since $\phi$ is an $\SL$-embedding, $M_n(\tilde K)$ and $M_n(\tilde L)$ induce the same $\SL$-structure $\tilde R$ on $R$ and $\phi|_R:\tilde R\lra \tilde K$, $\psi|_R:\tilde R\lra \tilde L$ are embeddings of $\SL$-structures. Since $\tilde F$ has quantifier elimination there are $\tilde \Omega\equiv \tilde F$ and $\SL$-embeddings $\rho:\tilde K\lra\tilde \Omega,\ \delta:\tilde L\lra\tilde \Omega$ such that $\rho\circ \phi|_R=\delta\circ \psi|_R$. We write ${\bar \rho},{\bar \delta}$ for the induced $\SL(\tr,\inv)$-embeddings as in \ref{AmalgTraces} and consider the diagram \begin{center} \begin{tikzcd}[row sep=3ex,column sep=1ex] & (M_n(\tilde \Omega),\tr_\Omega,X\mapsto X^t) \\ {} \\ (M_n(\tilde K),\tr_K,X\mapsto X^t) \arrow[uur, "\bar \rho"] && (M_n(\tilde L),\tr_L,X\mapsto X^t) \arrow[uul, "\bar \delta"'] \\ & \SU \arrow[ul,"\phi"'] \arrow[ur,"\psi"] \\ \\ &(R,f|_R) \arrow[uuul,"\phi|_R"] \arrow[uuur,"\psi|_R"'] \arrow[uu, hook] \end{tikzcd} \end{center} Notice that in general only the outer square in this diagram commutes. Since $\SU$ is a finitely generated $\SL$-structure, there are $X_1,\ldots,X_d\in U$ such that $U$ is the ring generated by $X_1,\ldots,X_d$. \par \noindent \claim There is some orthogonal matrix $O\in M_n(\Omega)$ such that for all $i\in\{1,\ldots,d\}$ we have \[ O\mal {\bar \rho}(\phi(X_i))\mal O^t={\bar \delta}(\psi(X_i)). \] \begin{proof} We write $Y_i={\bar \rho}(\phi(X_i))$ and $Z_i={\bar \delta}(\psi(X_i))$. To see the claim we use (2), by which it suffices to show that for every word $w$ in $x_1,\ldots,x_d,x_1^t,\ldots,x_d^t$ we have \[ \tr_\Omega(w(Y_1,\ldots,Y_d,Y_1^t,\ldots,Y_d^t))=\tr_\Omega(w(Z_1,\ldots,Z_d,Z_1^t,\ldots,Z_d^t)). \] Let $X=w(X_1,\ldots,X_d,h(X_1),\ldots,h(X_d))\in U$ (the degree bound $D$ is used to transfer (2) from $\tilde F$ to $\tilde \Omega$). By \ref{AmalgTraces}\ref{it:amal2} we know that \[ \tr_\Omega({\bar \rho}(\phi(X)))=\tr_\Omega({\bar \delta}(\psi(X))). \] Since ${\bar \rho}$ and $\phi$ respect the function symbol for the adjoint we see that \begin{align*} {\bar \rho}(\phi(X))&={\bar \rho}(\phi(w(X_1,\ldots,X_d,h(X_1),\ldots,h(X_d))))\cr &=w({\bar \rho}(\phi (X_1)),\ldots,{\bar \rho}(\phi (X_d)),{\bar \rho}(\phi (X_1))^t,\ldots,{\bar \rho}(\phi (X_d))^t)\cr &=w(Y_1,\ldots,Y_d,Y_1^t,\ldots,Y_d^t). \end{align*} Similarly, ${\bar \delta}(\psi(X))=w(Z_1,\ldots,Z_d,Z_1^t,\ldots,Z_d^t)$, establishing the claim. \end{proof} \noindent Now take an orthogonal $O\in M_n(\Omega)$ as in the claim and observe that the map $\gamma:M_n(\Omega)\lra M_n(\Omega),$ $X\mapsto O\mal X\mal O^t$ preserves traces, adjoints of matrices and the $\SL$-structure of $M_n(\tilde \Omega)$. Hence $\gamma$ is an $\SL(\tr,\inv)$-automorphism of $(M_n(\tilde \Omega),$ $\tr_\Omega,X\mapsto X^t)$. \par Consequently, by the claim, $\gamma\circ {\bar \rho}\circ \phi={\bar \delta}\circ \psi$. This shows that the maps $\gamma\circ {\bar \rho}\circ\phi $ and ${\bar \delta}\circ \psi$ form an amalgamation of the $\SL(\tr,\inv)$-structures $\SM$ and $\SN$ over the $\SL(\tr,\inv)$-embeddings $\phi$ and $\psi$. \par \smallskip\noindent (1)$\Ra $(3) is a weakening, see \ref{QEisMCplusAP}. \par \smallskip\noindent (3)$\Ra $(2) By a standard compactness argument it suffices to show that (2) holds without the degree bound for all $\tilde K\equiv \tilde F$. \par Let $\SU$ be the $\SL(\tr,\inv)$-substructure of $M_n(\tilde K)$ generated by $K\mal I_n$ and the $X_i$. Let $U$ be the ring underlying $\SU$. Hence $U$ is generated as a $K$-algebra by all words in the $X_i,X_i^t$. Let $\phi:U\to M_n(K)$ be the identity mapping and let $\psi:U\to M_n(K)$ be the $K$-algebra homomorphism that maps $X_i$ to $Y_i$ and $X_i^t$ to $Y_i^t$. \par We claim that $\psi$ is an $\SL(\tr,\inv)$-homomorphism. Firstly, $\psi$ is well defined: It suffices to show that for every noncommutative polynomial $p(x,x^t)$ with coefficients in $K$ and $p(X,X^t)=0$, we have $p(Y,Y^t)=0$. By \ref{lem:CCt} we know $\tr(p(X,X^t)^tp(X,X^t))=0$. But the left-hand side of this equation is simply a linear combination of traces of words in the $X,X^t$. Hence by the assumption on traces, $\tr(p(Y,Y^t)^tp(Y,Y^t))=0$. Thus $p(Y,Y^t)=0$ by \ref{lem:CCt} again. It is clear that $\psi$ is an $\SL(\tr,\inv)$-embedding. \par \smallskip Now we amalgamate as asserted in (3). There are an elementary extension $\tilde \Omega$ of $\tilde K$ and an $\SL$-embedding $\bar\rho:M_n(\tilde K)\to M_n(\tilde \Omega)$, preserving $\tr$ and $\inv$ such that $\psi(u)=\bar\rho(u)$ for all $u\in \SU$. Since $\bar\rho$ is compatible with the traces it is a $K$-algebra homomorphisms. Hence by the Skolem-Noether theorem (see \cite[Thm 4.46]{Bresar2014}), there is some invertible $Z\in M_n(\Omega)$ with \[ \bar\rho(X)=Z^{-1}XZ\quad \text{for all }X\in M_n(K). \] Now, \[ Z^{-1}X^tZ=\bar\rho(X^t)=(\bar\rho(X))^t = (Z^{-1}XZ)^t = Z^tX^t(Z^{-1})^t=Z^tX^t(Z^t)^{-1}, \] whence $ZZ^tX^t=X^tZZ^t$ for all $X$. Hence $ZZ^t$ is central and there is some $\lambda\in\sum \Omega^2$ with $ZZ^t=Z^tZ=\lambda I_n$. \par By the commutativity in the amalgamation diagram we know $Z^{-1} X_i Z = Y_i$ for all $i$. Since $\Omega$ is Pythagorean we also know that $\lambda$ is a square and so we may replace $Z$ by $\frac{Z}{\sqrt\lambda}$ and assume $\lambda=1$. Hence $O=Z^{-1}$ is an orthogonal matrix with coefficients in $\Omega$ satisfying $O^t X_iO=Y_i$ for all $i$. Since $\Omega$ is an elementary extension of $K$ we may find such an $O$ with coefficients in $K$ as well. \end{proof} \noindent We next identify a large class of fields with the Specht property, namely fields that can be written as intersections of real closed fields. We refer to \cite{Craven1980} for a systematic study of such fields. In \cite{MSV1993} the authors say such fields satisfy the principal axis property: every symmetric matrix over $F$ is orthogonally similar to a diagonal matrix over $F$. Notice that all fields that can be written as intersections of real closed fields are Pythagorean and by \cite[III, \S 1, Thm. 1]{Becker1978}, every hereditarily Pythagorean field is the intersection of real closed fields. \begin{FACT}{Proposition}{prop:Specht} Suppose the field $F$ is an intersection of real closed fields. Then $F$ has the Specht property for transposition. \par More precisely, given two $d$-tuples of $n\times n$ matrices $X,Y\in M_n(F)^d$ with \[ \tr(w(X,X^t))=\tr(w(Y,Y^t)) \] for all words $w$ in $x,x^t$ of length $\leq n^2$, there is some $O\in M_n(F)$ with $OO^t=I_n$ and $O^tX_iO=Y_i$ for all $i$. \end{FACT} \begin{proof} By \ref{Specht}, for every real closed field $R\supseteq F$ there is an orthogonal matrix $U\in M_n(R)$ with $U^tX_iU=Y_i$. \par Consider the system of linear equations $PX_i=Y_iP$ and $PX_i^t=Y_i^tP$ for $i=1,\ldots,d$. It has solutions $P$ with nonzero determinant in every real closed field extension of $F$ by the above, so it must have a solution $P\in M_d(F)$ that is invertible. Hence $P^{-1}X_iP=Y_i$ and $P^{-1}X_i^tP=Y_i^t$ for all $i$. In particular, \[ P^{-1}X_i^tP=Y_i^t=\left(P^{-1}X_iP\right)^t = P^t X_i^t (P^t)^{-1}, \] whence $PP^t$ commutes with all $X_i$ and $X_i^t$. \par Since $F$ has the principal axis property, we can diagonalize $PP^t$. There is an orthogonal matrix $V\in M_n(F)$ and a diagonal matrix $D\in M_n(F)$ with $V^t PP^t V = D$. By construction, each entry of $D$ is a (sum of) square(s). We thus find a diagonal matrix $\sqrt D\in M_n(F)$ with $\sqrt D^2=D$. Let $H:=V\sqrt DV^t\in M_n(F)$. Then \[ H^2=V\sqrt DV^tV\sqrt DV^t= V \sqrt D^2 V^t=VDV^t=PP^t, \] i.e., $H$ is the symmetric square root of $PP^t$. Thus by standard linear algebra, it commutes with all elements that commute with $PP^t$. \par Set $O=H^{-1}P$. Then \[ O^tO= P^tH^{-1}H^{-1}P = P^t H^{-2} P = P^t (PP^t)^{-1} P = P^t P^{-t} P^{-1} P = I, \] so $O\in M_n(F)$ is an orthogonal matrix. Further, \[ O^t X_i O = O^{-1} X_i O = P^{-1} H X_i H^{-1} P = P^{-1} X_i P = Y_i, \] as desired. \end{proof} \begin{FACT}{Corollary}{cor:PAP} Let $F$ be an intersection of real closed fields and let $\tilde F$ be an expansion of $F$ in a language $\SL$ extending the language of rings. Suppose $\tilde F$ has quantifier elimination in $\SL$. Let $\SL(\tr,\inv)$ be the extension of $\SL$ by two new unary function symbols. Then the structure $(M_n(\tilde F),\tr_F, X\mapsto X^t)$ has quantifier elimination in $\SL(\tr,\inv)$. \end{FACT} \begin{proof} Immediate from Theorem \ref{MnTrHqe} and Proposition \ref{prop:Specht}. \end{proof} \begin{FACT}{void}{Sylvester}\textbf{An application: Sylvester's equation} A famous matrix equation from control theory is Sylvester's equation \cite{BR97}, $AX-XB=C$ for some $n\in\N$ and $n\times n$ real (or complex) matrices $A,B,C$. By the Sylvester-Rosenblum theorem, given $A,B$ there is a unique solution $X$ for every $C$ iff the spectra of $A$ and $B$ are disjoint, and by the quantifier elimination proved in Corollary \ref{cor:PAP} (or \ref{ConjPrC} below, for the complex case), this can be expressed quantifier free in $A,B$ purely in terms of the trace and (conjugate) transpose. \end{FACT} \subsection{Trace and transposition are needed} \labelx{SectionExamples} We present three examples establishing the optimality of \ref{MnTrHqe}. The first example shows that we cannot omit the trace. \begin{FACT}{Example}{lem:noTrace} Let $K$ be a field of characteristic zero. Let \[ X_1=\begin{pmatrix} 1 \\ & 2 \\ &&2 \end{pmatrix},\qquad X_2=\begin{pmatrix} 1 \\ & 1 \\ &&2 \end{pmatrix}\in M_3(K). \] Let $U$ be the unital subring of $M_3(\Z)$ generated by $X_1$. Consider the ring homomorphism $\psi:U\to M_3(K)$ defined by \[ X_1\mapsto X_2, \] and let $\phi:U\to M_3(K)$ be the inclusion mapping. Then the following diagram cannot be amalgamated: \par \hfil$ \begin{tikzcd} M_3(K) & & M_3(K) \\ & U \arrow[hookrightarrow,lu,"\phi"] \arrow[ru,"\psi"'] \end{tikzcd} $ \par \noindent (Notice that $\phi$ and $\psi$ also respect the transposition, since all $X \in U$ are symmetric.) \end{FACT} \begin{proof} Notice that it suffices to verify the claim for $K=L=\Q$. Firstly, the map $\psi$ is well-defined, since the minimal polynomial of $X_1$ is $(t-1)(t-2)$ and is equal to the minimal polynomial of $X_2$. Now assume $M_3(\Omega)$ is an amalgamation of $\psi $ and $\varphi$ over $U$, and the following diagram commutes: \par \hfil \begin{tikzcd} & M_3(\Omega)\\ M_3(\Q) \arrow[ur,"\bar\rho"] & & M_3(\Q) \arrow [ul,"\bar\delta"'] \\ & U \arrow[hookrightarrow,lu,"\phi"] \arrow[ru,"\psi"'] \end{tikzcd} \par \noindent Then $\bar\rho,{\bar \delta}:M_3(\Q)\to M_3(\Omega)$ are ring homomorphisms. By the Skolem-Noether theorem (see, e.g. \cite[Theorem 4.46]{Bresar2014}), there are invertible matrices $W,V\in M_3(\Omega)$ such that \[\bar\rho(Y)= W^{-1} YW, \quad \bar\delta(Y)=V^{-1}YV\] for all $Y\in M_3(\Q)$. Then \[ V^{-1}X_2V= \bar\delta(X_2)=\bar\delta(\psi(X_1)) =\bar\rho(\varphi(X_1))= \bar\rho(X_1)=W^{-1}X_1V, \] yielding \[ X_2= (VW^{-1})\ X_1\ (VW^{-1})^{-1}. \] However, this is not possible because $X_1$ and $X_2$ are not similar; for example they have different characteristic polynomials. \end{proof} \par \noindent The second example shows that we cannot omit transposition in \ref{MnTrHqe}. \begin{FACT}{Example}{} Let $K$ be a field of characteristic zero. Consider \[ X_1=\begin{pmatrix} 0 & 1 & 0 &0 \\ 0 & 0 & 0 &0 \\ 0 & 0 & 0 &1 \\ 0 & 0 & 0 &0 \end{pmatrix},\qquad X_2=\begin{pmatrix} 0 & 1 & 0 &0 \\ 0 & 0 & 0 &0 \\ 0 & 0 & 0 &0 \\ 0 & 0 & 0 &0 \end{pmatrix}\in M_4(K). \] Let  $U$ be the unital subring of $M_4(\Z)$ generated by $X_1$. Consider the ring homomorphism $\psi:U\to M_4(K)$ defined by \[ X_1\mapsto X_2, \] and let $\phi:U\to M_4(K)$ be the identity mapping. Then $\psi$ and $\phi$ cannot be amalgamated over $U$. (Notice that $\phi$ and $\psi$ also respect the trace functions.) \end{FACT} \begin{proof} Again, it suffices to verify the claim for $K=L=\Q$. Note that $\psi$ is well-defined since the minimal polynomial of $X_1$ and of $X_2$ is $t^2$. Now assume $M_4(\Omega)$ amalgamates $\phi$ and $\psi$ over $U$. As in \ref{lem:noTrace} this leads to $X_1$ being conjugate to $X_2$ (over $\Omega$ and thus over $\Q$). However, this is impossible since $X_1$ and $X_2$ are not similar; for example $\dim\ker (X_1)=2\neq3=\dim\ker(X_2)$. \end{proof} \par \noindent By \ref{cor:PAP}, the structure $(M_n(\R),\leq ,\tr_\R, X\mapsto X^t)$ has quantifier elimination in $\SL_\mathrm{ri}(\leq,\tr,\inv)$. The third example shows that $(M_n(\C),\tr_\C, X\mapsto X^t)$ does not have quantifier elimination in $\SL_\mathrm{ri}(\tr,\inv)$. \begin{FACT}{Example}{} Complex matrices with the trace and transpose do not admit quantifier elimination. For the same reasons as above it suffices to show there exist symmetric order two nilpotents with different rank. For this we take $N_1$ to be the rank one outer product $N_1=uu^t$ with $u=\begin{pmatrix}1&i&0&0\end{pmatrix}^t$ and we let $N_2$ be the symmetric order two nilpotent \[ N_2=\left(\begin{array}{rrrr} 0&1&0&-i\\ 1&0&-i&0\\ 0&-i&0&-1\\ -i&0&-1&0. \end{array}\right). \] \end{FACT} \subsection{The simultaneous conjugacy problem} \labelx{QEoverC} \numberwithin{theorem}{subsection}\ \begin{FACT}{void}{} As in the proof of \ref{MnTrHqe}(2)$\Ra $(1), using the complex Specht property (see \ref{Specht}), one can establish that the theory of $(M_n(\C),\leq,\tr_{\C}, X\mapsto X^*)$ has quantifier elimination; here $\leq$ is the order on the symmetric center $\R\cdot I_n$. The underlying expansion of the field $\C$ is $\tilde \C:=(\C,z\mapsto \overline{z},\leq)$, where $\leq$ is the order on $\R$ and $\overline{z}$ is complex conjugation. Since $\R$ is not definable in the field $\C$, the structure $\tilde \C$ is a proper expansion of $\C$. Conversely, the field $\R$ obviously defines the structure $\tilde \C$; hence the complex version of \ref{MnTrHqe} is a statement about the real field. \end{FACT} \begin{FACT}{void}{ConjPrC} The question on whether a natural \textit{definable} expansion of the ring $M_n(\C)$\footnote{Hence an expansion of $M_n(\C)$, which interprets the new symbols by sets and functions that are definable in $M_n(\C)$.} has quantifier elimination is tightly related to a ``hopeless'' open problem in invariant theory \cite{LeBruyn1995,lBrPro1987,GelPom1969}. Namely the classification of $d$-tuples of $n\times n$ matrices under simultaneous conjugation by $\GL_n(\C)$, i.e., understanding the quotient $M_n(\C)^d/\GL_n(\C)$. Alternately, in algebraic language, one is interested in a canonical form for tuples of matrices under simultaneous conjugation, a role played by the Jordan canonical form in the case $d=1$. A relaxation of the problem asks for a set of invariants that separate the orbits. \par In model theoretic terms this can be phrased as follows. Let $M$ be the ring $M_n(\C)$ and fix $d\in\N$. We write $\sim_d$ for the simultaneous similarity relation on $M^d$, i.e., $X\sim_dY$ if and only if there is $Z\in \GL_n(\C)$ with $X=Z^{-1}YZ$. Then $\sim_d$ is a 0-definable equivalence relation and by elimination of imaginaries of the field $\C$ (cf. \cite[Thm. 4.4.6]{Hodges1993}), there is a 0-definable function $f_d:M^d\lra M^k$ for some $k$ such that ${\bar X}\sim{\bar Y}\iff f_d({\bar X})=f_d({\bar Y})$. If we add names for all the $f_d$ to the language of rings, one can prove quantifier elimination of the resulting expansion of $M$ just like in the proof of \ref{MnTrHqe}(2)$\Ra $(1); the sequence of the $f_d$ substitutes the role of the transposition and the trace {\color{black}(most crucially the claim in that proof becomes just the defining property of the $f_d$'s)}. \par {\color{black}A caveat here is that the functions $f_d$ are not explicit.} In \cite{Friedlan1983} functions $f_d$ as above are explicitly constructed, up to a finite number of exceptions. Alternatively one can use techniques from Gr\"obner bases to construct them explicitly (without exceptions). This is work in progress and will be published in another paper. \end{FACT} \smallskip We would also like to point out that after submission of this paper, a subsequent paper \cite{DKMV2023} -- published in 2023 -- has identified a concrete set of separating invariants for $M^d$. Namely, a $d$-tuple of $n\times n$ matrices $(A_1,\ldots,A_d)$ is up to simultaneous similarity uniquely determined by ranks of linear matrix pencils \[ \mathrm{rank}\, (I_n \otimes T_0 + A_1 \otimes T_1 + \cdots + A_d \otimes T_d), \] where the $T_j$ run through $m\times m$ matrices with $m\leq dn$. \section{Undecidability of dimension-free matrices} \labelx{SectionUndecidable} \numberwithin{theorem}{subsection} \par \labelx{SectionDimFree} \noindent We now turn to model theoretic properties of dimension-free matrices. In Section \ref{sectionAppStructures} we present six natural algebraic structures capturing the set of all matrices of all sizes over a given field and prove that all of them are undecidable. This is based on undecidability of finite groups, which is reviewed in Section \ref{sectionAppGroups}, suitable for our purpose. As a general reference for elementary properties of classes of finite groups in relation to decidability questions, we refer to \cite[Section 6.3]{BunMik2004}. \subsection{The universal Horn theory of finite groups}\labelx{sectionAppGroups} \noindent Throughout, $\SL_\mathrm{gr}$ denotes the language $\{\cdot, \ ^{-1},e\}$ of groups and $T_\mathrm{fin}$ denotes the $\SL_\mathrm{gr}$-theory of finite groups. Hence \[ T_\mathrm{fin}=\{\phi \st \phi\text{ an } \SL_\mathrm{gr}\text{-sentence with }G\models \phi\text{ for every finite group }G\}. \] Further, $T_{\mathrm{fin},\forall}$ denotes the \notion{universal theory of finite groups}, hence all sentences in $T_\mathrm{fin}$ of the form \[ \forall x_1,\ldots,x_n\bigwedge_{\lambda=1}^r \biggl(\bigwedge _{j=1}^ms_{\lambda j}=e\ \lra \ \bigvee _{i=1}^k t_{\lambda i}=e\biggr), \] where $r,m,k\in \N_0$, $r\geq 1$ and $s_{\lambda j},t_{\lambda i}$ are $\SL_\mathrm{gr}$-terms in the free variables $x_1,\ldots,x_n$ (aka ``words in the $x_i$ and $x_i^{-1}$''). A \notion{universal Horn sentence} of $\SL_\mathrm{gr}$ is a sentence of the form \[ \forall x_1,\ldots,x_n\biggl(\bigwedge _{j=1}^ms_{j}=e\ \lra \ t=e\biggr), \] where $m\in \N_0$ and $s_{j},t$ are $\SL_\mathrm{gr}$-terms. We write $T_{\mathrm{fin},\mathrm{H-}\forall}$ for the set of all universal Horn sentences in $T_{\mathrm{fin},\forall}$ and call it the \notion{universal Horn theory of finite groups}. \par Notice that by the shape of the sentences in $T_{\mathrm{fin},\forall}$ and in $T_{\mathrm{fin},\mathrm{H-}\forall}$, every subgroup of a model of $T_{\mathrm{fin},\forall}$, $T_{\mathrm{fin},\mathrm{H-}\forall}$ is again a model of $T_{\mathrm{fin},\forall}$, $T_{\mathrm{fin},\mathrm{H-}\forall}$ respectively. \begin{FACT}{Fact}{Slobodskoi} (cf. \cite{Slobod1981})\\ The universal Horn theory of finite groups is undecidable. More precisely: $T_{\mathrm{fin},\mathrm{H-}\forall}$ is not a recursive subset of the set of $\SL_\mathrm{gr}$-sentences. The same is then obviously true for $T_{\mathrm{fin},\forall}$. \end{FACT} \begin{FACT}{Definition}{defnSatiated} We call a class $\CK$ of groups \notion{satiated} if \begin{enumerate}[(a)] \item Every finite group embeds into some member of $\CK$, and, \item Every member of $\CK$ is a model of the universal Horn theory of finite groups. \end{enumerate} \noindent Let $\CR$ be any first order structure in an arbitrary language $\SL$. We call $\CR$ \notion{satiated} if $\CR$ has a uniform interpretation of a satiated set of groups. This means that there are $k,n\in\N$ and an $\SL$-formula $\mu({\bar x_1},{\bar x_2},{\bar y},{\bar z})$, where ${\bar x_1},{\bar x_2},{\bar y}$ are $n$-tuples and ${\bar z}$ is a $k$-tuple such that \begin{enumerate}[(a)] \item for every ${\bar a}\in\CR^k$, the subset defined by $\mu({\bar x_1},{\bar x_2},{\bar y},{\bar a})$ in $\CR^{3n}$ is the graph of multiplication of a group $G_{\bar a}$ with universe contained in $\CR^n$, and, \item the set of groups $\{G_{\bar a}\st {\bar a}\in \CR^k\}$ is satiated. \end{enumerate} \end{FACT} \begin{FACT}{Proposition}{SatiatedUndecidable} Any satiated structure is undecidable. \end{FACT} \begin{proof} The definition readily implies that the universal Horn theory of every satiated class $\CK$ (thus, all universal Horn $\SL_\mathrm{gr}$-sentences that are true in all $G\in \CK$) is the universal Horn theory of finite groups. Now suppose that $\CR$ is a decidable satiated structure. Take a formula $\mu$ as in \ref{defnSatiated}. It is then clear that there is a map $\phi\mapsto \tilde \phi$ from universal Horn sentences in $\SL_\mathrm{gr}$ to the set of $\SL$-sentences with recursive image such that $\phi\in T_{\mathrm{fin},\mathrm{H-}\forall}$ if and only if $\tilde \phi$ is true in $\CR$. But then $T_{\mathrm{fin},\mathrm{H-}\forall}$ is recursive, in contradiction to \ref{Slobodskoi}.\end{proof} \par \medskip\noindent Recall that a \textbf{linear group} is a group that can be embedded into some $\Gl_n(F)$ for some field $F$. \begin{FACT}{Theorem}{LinearGroupsAreUniversallyPseudofinite} Every linear group is a model of the universal theory of finite groups.\footnote{Note that each finite group $G$ embeds into $M_{|G|}(F)$ via the left regular representation $\lambda$ in such a way that $\tr(\lambda(g))=0$ for $g\neq1$ and $\lambda(g^{-1})=\lambda(g)^t$.} \par \end{FACT} \begin{proof} It suffices to show the claim for the group $G=\Gl_n(F)$ when $F$ is an algebraically closed field. If $F$ has characteristic $p>0$, then by completeness of the theory of algebraically closed fields of fixed characteristic we may assume that $F$ is the algebraic closure $\overline{\F_p}$ of $\F_p$. But then $G$ is the union of all the $\Gl_n(K)$, where $K$ runs through the finite fields of characteristic $p$. Since universal sentences are preserved by unions we get the assertion. When $F$ is of characteristic 0, then using {\L}o\'s's theorem, $G$ is elementarily equivalent to any non-principal ultraproduct of the $\Gl_n(\overline{\F_p})$, $p$ prime. Hence the result follows. \end{proof} \begin{FACT}{Corollary}{GlNSatiated} Let $\CK$ be any class of linear groups such that every finite group embeds into some member of $\CK$. Then $\CK$ is satiated. This, for example, is the case for any class of linear groups containing all the $\Gl_n(F)$ for some fixed field $F$. \end{FACT} \begin{proof} Immediate from \ref{LinearGroupsAreUniversallyPseudofinite}.\end{proof} \subsection{Applications to dimension-free matrices}\labelx{sectionAppStructures} There are various ways how the collection of all square matrices of arbitrary (finite) size over a field can be given an algebraic structure. We present six such constructions and show that each of them is undecidable. In the realm of infinite matrix theory in the sense of Poincar\'e (cf. \cite{Bernko1968} and \cite{Cooke1950}), one can find many constructions containing all finite square matrices. But then either one does not have a handle on the finitely sized matrices, or one of the constructions below will be interpretable. \begin{FACT}{:Dimension-free matrices with partial operations.}{partialOp} Let $F$ be a field and let $\CR_1,\CR_2$ be the following structures in a language $\SL=\{R\}$ for a ternary relation symbol $R$. The universe of $\CR_1$ is the disjoint union of all the $\Gl_n(F)$. The relation symbol $R$ is interpreted in $\CR_1$ as the union of the graphs of all the multiplication maps $\Gl_n(F)\times \Gl_n(F)\lra \Gl_n(F)$. The relation symbol $R$ is interpreted in $\CR_2$ as the union of the graphs of all the multiplication maps $M_n(F)\times M_n(F)\lra M_n(F)$. \par Then $\CR_1,\CR_2$ are satiated, hence undecidable by \ref{SatiatedUndecidable}. The formula $\mu$ that uniformly interprets the satiated set $\{\Gl_n(F)\st n\in \N\}$ in $\CR_1$ is the formula \[ \exists u\, R(x_1,z,u)\land \exists u\, R(x_2,z,u)\land R(x_1,x_2,y). \] For $\CR_2$ we take the formula $\mu(x_1,x_2,y,z)\ \&\ ``x_1,x_2\text{ are invertible''}$, where ``$x$ invertible'' stands for the $\SL$-formula expressing that $x$ is invertible in the semigroup of all $u$ for which $u\mal x$ is defined. \end{FACT} \begin{FACT}{Lemma}{SubSemiIsLinear} Let $F$ be a field and let $S$ be a subsemigroup of $M_n(F)$. If $S$ is a group, then $S$ is isomorphic to a subgroup of $\Gl_m(F)$ for some $m\leq n$. In particular, $S$ is a linear group. \end{FACT} \begin{proof} Let $I$ be the neutral element of $S$. Then $I$ is idempotent and there is some $P\in \Gl_n(F)$ such that $P^{-1}\mal I\mal P$ is of the form \[ E'= \begin{pmatrix} E & 0 \\ 0 & 0 \end{pmatrix}, \] \noindent where $E$ is the identity matrix of $M_m(F)$ for some $m\leq n$. Let $\sigma:M_n(F)\lra M_n(F);\ \sigma (X)=P^{-1}\mal X\mal P$. Then $\sigma $ is an automorphism of $M_n(F)$ and as $I\mal X\mal I=X$ we get $E'\mal \sigma (X)\mal E'=\sigma (X)$ for all $X\in S$. However, matrices with this property are all of the form \[ Y'= \begin{pmatrix} Y & 0 \\ 0 & 0 \end{pmatrix}, \] for some $Y\in M_m(F)$. If we embed $M_m(F)$ into $M_n(F)$ by mapping $Y$ to $Y'$, we see that $\sigma $ maps $S$ into $M_m(F)$. Hence $S$ is isomorphic to a subgroup of $\Gl_m(F)$. \end{proof} \begin{FACT}{+Finite rank infinite matrices.}{} Let $F$ be a field and let $\CR$ be the semigroup of all $\N\times \N$-matrices with finite support and multiplication as operation. Then $\CR$ is a satiated structure and is thus undecidable by \ref{SatiatedUndecidable}. \end{FACT} \begin{proof} We consider $M_n(F)$ as the subsemigroup of $\CR$ consisting of all $n\times n$-matrices sitting in the corner of $\CR$. We give a uniform definition of a satiated class of linear groups in $\CR$ using a formula $\mu$ in the language $\{\cdot \}$ of semigroups, as explained in \ref{defnSatiated}. For $X\in \CR$, consider the set \[ \CC(X)=\{Y\in \CR\st \forall Z\in \CR\, \bigl((X\mal Z=0\ra Y\mal Z=0)\ \&\ (Z\mal X=0\ra Z\mal Y=0)\bigr)\}. \] It is easy to see that $\CC(X)\subseteq M_n(F)$ for $X\in M_n(F)$ and that $\CC(X)=M_n(F)$ for $X\in \Gl_n(F)$. \par Let $\psi(z_1,z_2)$ be an $\{\cdot\}$-formula such that $\psi$ holds at $(X,I)\in \CR^2$ in $\CR$ just if the set \begin{quote} $\CG(X,I)=\{Y\in \CC(X)\st \exists Z\in \CC(X)\ Y\mal Z=Z\mal Y=I\}$ \end{quote} is a group with neutral element $I$. Then the formula $\phi(x,z_1,z_2)$ defined as \[ (\psi(z_1,z_2)\ra x\in \CG(z_1,z_2))\ \&\ (\lnot \psi(z_1,z_2)\ra x=0) \] has the following properties for all $(X,I)\in \CR^2$: \begin{enumerate}[(a)] \item The set of all $Y\in\CR$ with $\CR\models\phi(Y,X,I)$ is a linear group (use \ref{SubSemiIsLinear}). \item If $X\in \GL_n(F)$ and $I=I_n$, then set of all $Y\in\CR$ with $\CR\models\phi(Y,X,I)$ is $\Gl_n(F)$. \end{enumerate} It is now standard to write down a $\{\cdot\}$-formula $\mu(x_1,x_2,y,z_1,z_2)$ that uniformly defines a satiated class of groups (also invoke \ref{GlNSatiated}). \end{proof} \begin{FACT}{+Products.}{} If $(G_i\st i\in I)$ is a satiated family of groups, then $\prod_{i\in I}G_i$ is undecidable, in fact the universal Horn theory of that product is undecidable. Hence by \ref{GlNSatiated}, for any field $F$, the group $\prod_{n\in \N}\Gl_n(F)$ is undecidable, and consequently so is the semigroup $\prod_{n\in \N}M_n(F)$ (observe that $\prod_{n\in \N}\Gl_n(F)$ is the set of invertible elements of $\prod_{n\in \N}M_n(F)$). \end{FACT} \begin{proof} We write $P=\prod_{i\in I}G_i$ and show that $P$ satisfies exactly the same universal Horn sentences as the those satisfied by all finite groups. Then \ref{Slobodskoi} gives the assertion. \par \smallskip\noindent As a product, $P$ satisfies all universal Horn sentences that are true in all $G_i$ and so $P$ satisfies all universal Horn sentences that are true in all finite groups. \par \smallskip\noindent Conversely, let $\phi$ be a quantifier-free Horn formula \[ \bigwedge _js_j=e\ra t=e \] in $l$ free variables and assume $P\models \forall \phi$. Let $H$ be a finite group and suppose $H\models \bigwedge _js_j(h_1,\ldots,h_l)=e$. Fix some $i_0\in I$ and an embedding $\iota :H\into G_{i_0}$. We define $X_1,\ldots,X_l\in P$ by \[ X_{j,i}= \begin{cases} \iota(h_j) & \text{if }i=i_0, \cr e & \text{if }i\neq i_0. \end{cases} \] It is clear that $G_i\models \bigwedge _js_j(X_{1,i},\ldots,X_{l,i})=e$ for all $i\in I$. Hence \[ P\models \bigwedge _js_j(X_{1},\ldots,X_{l})=e \] and so $P\models t(X_{1},\ldots,X_{l})=e$. Looking at the $i_0^\mathrm{th}$ component we see that $H\models t(h_1,\ldots,h_l)=e$ as required. \end{proof} \begin{FACT}{+Ultraproducts.}{} For any field $F$ and any non-principal ultrafilter $\DU$ on $\N$, the universal Horn theory of the ultraproduct $\prod_{n\in \N}\Gl_n(F)/\DU$ is the universal Horn theory of finite groups, and is thus undecidable. Since the natural map \[ \prod _n\Gl_n(F)/\DU\lra (\prod _nM_n(F)/\DU)^\times \] is an isomorphism, the semigroup $\prod _nM_n(F)/\DU$ is undecidable as well. \end{FACT} \begin{proof} \smallskip Let $G_\infty=\prod _n\Gl_n(F)/\DU$, for some non-principal ultrafilter $\DU$. If $\phi $ is a universal sentence, true in all finite groups, then by \ref{LinearGroupsAreUniversallyPseudofinite} it is true in all $\Gl_n(F)$ and so it is also true in $G_\infty$. \par Conversely if $G_\infty\models \phi$, then $\phi $ is true in all finite groups: Let $H$ be a finite group and let $N\in \N$ be such that $\Gl_n(F)$ contains an isomorphic copy $H_n$ of $H$ for all $n\geq N$. Since $\Gl_n(F)\models \phi $ for arbitrarily large $n$ and $\phi $ is universal, $\phi $ is also true in $H_n$. \par Hence the universal theory of the ultraproduct is $T_{\mathrm{fin},\forall}$. Now use \ref{Slobodskoi}. \end{proof} \begin{FACT}{+Direct Limits.}{DirLim} Let $F$ be a field. For $n\in\N$ let $f_{n}:M_{2^n}(F)\lra M_{2^{n+1}}(F)$ be the ring homomorphism that sends $X$ to $\begin{pmatrix} X & 0 \\ 0 & X \end{pmatrix}$. We consider the direct limit $\varinjlim M_{2^n}(F)$ induced by the $f_n$. \par Then for every infinite field $F$, the ring $\varinjlim M_{2^n}(F)$ is undecidable. In fact, it interprets the weak monadic second order logic of $F$. \end{FACT} \begin{proof} By the weak monadic second order logic of the field $F$ we mean the following first order structure $W$ expanding the poset $P$ of finite subsets of $F$: We identify $F$ with the subset $\{\{a\}\st a\in F\}$ of $P$ and expand $P$ by the graph of addition and multiplication of $F$; for details see, for example, \cite{Bauval1985} or \cite[Section 2]{Tressl2017}. \par We now show that $W$ is interpretable in $\varinjlim M_{2^n}(F)$. Firstly, we identify $F$ with the center of $\varinjlim M_{2^n}(F)$, which is 0-definable therein. If $X\in \varinjlim M_{2^n}(F)$, then let $\sigma (X)$ be the set of all central elements $\Lambda\in \varinjlim M_{2^n}(F)$ such that there is no $Y\in \varinjlim M_{2^n}(F)$ with $(X-\Lambda)\mal Y=I$. Hence $\sigma(X)$ is the finite set of eigenvalues of $X$. The map $\sigma$ is obviously 0-definable in $\varinjlim M_{2^n}(F)$. Further, if $X,Y\in \varinjlim M_{2^n}(F)$, then the property $\sigma(X)\subseteq \sigma(Y)$ is 0-definable in $\varinjlim M_{2^n}(F)$. \par The universe of $W$ then is the image of $\sigma$, i.e., the set $P$ of finite subsets of $F$ and the partial order on $P$ is interpretable in $\varinjlim M_{2^n}(F)$. On central elements, the map $\sigma$ is injective, hence the graph of addition and multiplication on the atoms of $W$ is interpretable in $\varinjlim M_{2^n}(F)$ as well. \par Hence $\varinjlim M_{2^n}(F)$ interprets $W$ and $W$ is well known to be undecidable, see for example \cite[2.5]{Tressl2017} for $\Char(F)=0$ and \cite[2.6]{Tressl2017} for $\Char(F)>0$. \end{proof} \begin{FACT}{+Row and column finite matrices.}{} Let $F$ be an infinite field and let $I$ be an infinite index set. Let $M_I(F)$ be the set of all $I\times I$ matrices $X$ such that all but a finite number of entries in each row and each column of $X$ are $0$. One checks that $M_I(F)$ is a ring under the ordinary definition of addition and multiplication. \par Then the ring $M_{I}(F)$ is undecidable. \end{FACT} \begin{proof} The interpretation used in the proof of \ref{DirLim} now gives the monadic second order theory of $F$, where second order quantifiers range over subsets of $F$ of size at most the cardinality of $I$. This is undecidable as well, see the proofs of \cite[2.5, 2.6]{Tressl2017}.\looseness=-1 \end{proof} \par \def\cprime{$'$}  \par 
\begin{thebibliography}{DKMV23} \bibitem[AM16]{AglMcC2016} J.~Agler and J.~E. McCarthy. \newblock The implicit function theorem and free algebraic sets. \newblock {\em Trans. Amer. Math. Soc.}, 368(5):3157--3175, 2016. \par \bibitem[Bau85]{Bauval1985} A.~Bauval. \newblock Polynomial rings and weak second-order logic. \newblock {\em J. Symbolic Logic}, 50(4):953--972 (1986), 1985. \par \bibitem[Bec78]{Becker1978} E.~Becker. \newblock {\em Hereditarily-{P}ythagorean fields and orderings of higher level}, volume~29 of {\em Monograf\'{i}as de Matem\'{a}tica [Mathematical Monographs]}. \newblock Instituto de Matem\'{a}tica Pura e Aplicada, Rio de Janeiro, 1978. \par \bibitem[Ber68]{Bernko1968} M.~Bernkopf. \newblock A history of infinite matrices. \newblock {\em Arch. History Exact Sci.}, 4(4):308--358, 1968. \newblock A study of denumerably infinite linear systems as the first step in the history of operators defined on function spaces. \par \bibitem[BM04]{BunMik2004} E.~I. Bunina and A.~V. Mikhalev. \newblock Elementary properties of linear groups and related problems. \newblock {\em J. Math. Sci. (N. Y.)}, 123(2):3921--3985, 2004. \newblock Algebra. \par \bibitem[BR97]{BR97} Rajendra Bhatia and Peter Rosenthal. \newblock How and why to solve the operator equation {$AX-XB=Y$}. \newblock {\em Bull. London Math. Soc.}, 29(1):1--21, 1997. \par \bibitem[Bre14]{Bresar2014} M.~Bre{\v{s}}ar. \newblock {\em Introduction to noncommutative algebra}. \newblock Universitext. Springer, Cham, 2014. \par \bibitem[BSZ09]{BaSaZu2009} R.~Bautista, L.~Salmer\'{o}n, and R.~Zuazua. \newblock {\em Differential tensor algebras and their module categories}, volume 362 of {\em London Mathematical Society Lecture Note Series}. \newblock Cambridge University Press, Cambridge, 2009. \par \bibitem[Coo50]{Cooke1950} R.~G. Cooke. \newblock {\em Infinite {M}atrices and {S}equence {S}paces}. \newblock Macmillan \& Co., Ltd., London, 1950. \par \bibitem[Cra80]{Craven1980} T.~C. Craven. \newblock Intersections of real closed fields. \newblock {\em Canad. J. Math.}, 32(2):431--440, 1980. \par \bibitem[DKMV23]{DKMV2023} Harm Derksen, Igor Klep, Visu Makam, and Jurij Vol\v{c}i\v{c}. \newblock Ranks of linear matrix pencils separate simultaneous similarity orbits. \newblock {\em Adv. Math.}, 415:Paper No. 108888, 20, 2023. \par \bibitem[DNT23]{DrNeTh2023} Tom Drescher, Tim Netzer, and Andreas Thom. \newblock On projections of free semialgebraic sets. \newblock {\em Adv. Geom.}, 23(2):207--214, 2023. \par \bibitem[FHS14]{FaHaSh2014} I.~Farah, B.~Hart, and D.~Sherman. \newblock Model theory of operator algebras {II}: model theory. \newblock {\em Israel J. Math.}, 201(1):477--505, 2014. \par \bibitem[Fri83]{Friedlan1983} S.~Friedland. \newblock Simultaneous similarity of matrices. \newblock {\em Adv. Math.}, 50(3):189--265, 1983. \par \bibitem[GfP69]{GelPom1969} I.~M. Gel\cprime~fand and V.~A. Ponomarev. \newblock Remarks on the classification of a pair of commuting linear transformations in a finite-dimensional space. \newblock {\em Funkcional. Anal. i Prilo\v{z}en.}, 3(4):81--82, 1969. \par \bibitem[HKM11]{HeKlMc2011} J.~W. Helton, I.~Klep, and S.~McCullough. \newblock Proper analytic free maps. \newblock {\em J. Funct. Anal.}, 260(5):1476--1490, 2011. \par \bibitem[Hod93]{Hodges1993} W.~Hodges. \newblock {\em Model theory}, volume~42 of {\em Encyclopedia of Mathematics and its Applications}. \newblock Cambridge University Press, Cambridge, 1993. \par \bibitem[KVV14]{KaVVin2014} D.~S. Kaliuzhnyi-Verbovetskyi and V.~Vinnikov. \newblock {\em Foundations of free noncommutative function theory}, volume 199 of {\em Mathematical Surveys and Monographs}. \newblock American Mathematical Society, Providence, RI, 2014. \par \bibitem[LB95]{LeBruyn1995} L.~Le~Bruyn. \newblock {\em Orbits of matrix tuples}. \newblock Universitaire Instelling Antwerpen. Department of Mathematics, 1995. \par \bibitem[LBP87]{lBrPro1987} L.~Le~Bruyn and C.~Procesi. \newblock {\'E}tale local structure of matrix invariants and concomitants. \newblock In {\em Algebraic groups Utrecht 1986}, pages 143--175. Springer, 1987. \par \bibitem[MSV93]{MSV1993} D.~Mornhinweg, B.~Shapiro, and K.~G. Valente. \newblock The principal axis theorem over arbitrary fields. \newblock {\em Amer. Math. Monthly}, 100(8):749--754, 1993. \par \bibitem[Pea62]{Pearcy1962} C.~Pearcy. \newblock A complete set of unitary invariants for operators generating finite {W}*-algebras of type {I}. \newblock {\em Pacific J. Math}, 12:1405--1416, 1962. \par \bibitem[Pro76]{Procesi1976} C.~Procesi. \newblock The invariant theory of {$n\times n$} matrices. \newblock {\em Adv. Math.}, 19(3):306--381, 1976. \par \bibitem[Pro07]{procesiBook} Claudio Procesi. \newblock {\em Lie groups}. \newblock Universitext. Springer, New York, 2007. \newblock An approach through invariants and representations. \par \bibitem[Put07]{Putinar2007} M.~Putinar. \newblock Undecidability in a free *-algebra. \newblock 2007. \par \bibitem[Rau10]{Rauten2010} W.~Rautenberg. \newblock {\em A concise introduction to mathematical logic}. \newblock Universitext. Springer, New York, third edition, 2010. \newblock With a foreword by Lev Beklemishev. \par \bibitem[Raz74]{Razmys1974} Ju.~P. Razmyslov. \newblock Identities with trace in full matrix algebras over a field of characteristic zero. \newblock {\em Izv. Akad. Nauk SSSR Ser. Mat.}, 38:723--756, 1974. \par \bibitem[Ros78]{Rose1978a} B.~I. Rose. \newblock Rings which admit elimination of quantifiers. \newblock {\em J. Symbolic Logic}, 43(1):92--112, 1978. \par \bibitem[Ros80]{Rose1980a} B.~I. Rose. \newblock On the model theory of finite-dimensional algebras. \newblock {\em Proc. London Math. Soc. (3)}, 40(1):21--39, 1980. \par \bibitem[Row80]{rowen} Louis~Halle Rowen. \newblock {\em Polynomial identities in ring theory}, volume~84 of {\em Pure and Applied Mathematics}. \newblock Academic Press, Inc. [Harcourt Brace Jovanovich, Publishers], New York-London, 1980. \par \bibitem[Sib68]{Sibir1968} K.~S. Sibirski\u\i. \newblock Algebraic invariants of a system of matrices. \newblock {\em Sibirsk. Mat. \v Z.}, 9:152--164, 1968. \par \bibitem[Slo81]{Slobod1981} A.~M. Slobodsko{\u\i}. \newblock Undecidability of the universal theory of finite groups. \newblock {\em Algebra i Logika}, 20(2):207--230, 251, 1981. \par \bibitem[Spe40]{Specht1940} W.~Specht. \newblock Zur {T}heorie der {M}atrizen. {II}. \newblock {\em Jber. Deutsch. Math. Verein.}, 50:19--23, 1940. \par \bibitem[Tre17]{Tressl2017} M.~Tressl. \newblock On the strength of some topological lattices. \newblock In {\em Ordered algebraic structures and related topics}, volume 697 of {\em Contemp. Math.}, pages 325--347. Amer. Math. Soc., Providence, RI, 2017. \par \bibitem[Voi10]{Voiculescu2010} D.-V. Voiculescu. \newblock Free analysis questions {II}: the {G}rassmannian completion and the series expansions at the origin. \newblock {\em J. Reine Angew. Math.}, 645:155--236, 2010. \par \bibitem[Wie62]{Wiegma1961} N.~A. Wiegmann. \newblock Necessary and sufficient conditions for unitary similarity. \newblock {\em J. Austral. Math. Soc.}, 2:122--126, 1961/1962. \par \end{thebibliography}
\end{document}